\documentclass[A4paper]{amsart}
\usepackage{amssymb, enumerate, xspace}
\usepackage[all]{xy}
\SelectTips{cm}{}

\numberwithin{equation}{section}

1

\newenvironment{enumeratea}
{\begin{enumerate}[\upshape (a)]}
{\end{enumerate}}

\newtheorem*{namedtheorem}{\theoremname}
\newcommand{\theoremname}{testing}

\newtheorem{theorem}{Theorem}[section]
\newtheorem{proposition}[theorem]{Proposition}
\newtheorem{proposition-definition}[theorem]
{Proposition-Definition}
\newtheorem{corollary}[theorem]{Corollary}
\newtheorem{lemma}[theorem]{Lemma}

\theoremstyle{definition}
\newtheorem{definition}[theorem]{Definition}

\newtheorem{examples}[theorem]{Examples}
\newtheorem{remark}[theorem]{Remark}

\theoremstyle{remark}

\newcommand\nome{testing}
\newcommand\call[1]{\label{#1}\renewcommand\nome{#1}}
\newcommand\itemref[1]{\item\label{\nome;#1}}
\newcommand\refall[2]{\ref{#1}~(\ref{#1;#2})}
\newcommand\refpart[2]{(\ref{#1;#2})}

\renewcommand\H{\operatorname{H}}

\newcommand\eqdef{\overset{\mathrm{def}} =}

\newcommand\into{\hookrightarrow}

\renewcommand\th{^\text{th}}

\def\displaytimes_#1{\mathrel{\mathop{\times}\limits_{#1}}}

\def\displayotimes_#1{\mathrel{\mathop{\bigotimes}\limits_{#1}}}

\newcommand\pic{\operatorname{Pic}}

\newcommand\spec{\operatorname{Spec}}

\newcommand{\GL}{\operatorname{GL}}

\newcommand{\PGL}{\operatorname{PGL}}

\newcommand{\SL}{\operatorname{SL}}

\newcommand{\gm}{\mathbb{G}_{\mathrm{m}}}

\newcommand\id{\operatorname{id}}

\newdir{ >}{{}*!/-5pt/@{>}}

\newcommand{\mmu}[1][r]{\boldsymbol{\mu}_{#1}}

\newcommand\cov[1][n,r,d]{\mathcal{H}(#1)}

\newcommand\covv[1][n,r,d]{\widetilde{\mathcal{H}}(#1)}

\newcommand\altcov[1][n,r,d]{\mathcal{H}'(#1)}

\newcommand\altaltcov[1][n,r,d]{\mathcal{H}''(#1)}

\newcommand\covsm[1][n,r,d]{\mathcal{H}_{\mathrm{sm}}(#1)}

\newcommand\covcomp[1][n,r,d]{\overline{\mathcal{H}}(#1)}

\newcommand\triplee[1][1,3;d_1, d_2]{\widetilde{\mathcal{H}}(#1)}

\newcommand\triple[1][1,3;d_1, d_2]{\mathcal{H}(#1)}

\newcommand\triplesm[1][1,3;d_1, d_2]{\mathcal{H}_{\mathrm{sm}}(#1)}

\newcommand\alttriple[1][1,3;d_1, d_2]{\mathcal{H}'(#1)}

\newcommand\covof[2][r]{\mathcal{H}(#2, #1)}

\newcommand\altcovof[2][r]{\mathcal{H}'(#2, #1)}

\newcommand\forms[1][n,rd]{\mathbb{A}(#1)}

\newcommand\pforms[1][n,rd]{\mathbb{P}(#1)}

\newcommand\formsnoto[1][n,rd]{\mathbb{A}_0(#1)}

\newcommand\formssm[1][n,rd]{\mathbb{A}_\mathrm{sm}(#1)}

\newcommand\tot{uniform cyclic cover\xspace}

\newcommand\tots{uniform cyclic covers\xspace}

\newcommand\Tots{Uniform cyclic covers\xspace}

\newcommand\Div{\operatorname{Div}}

\newcommand\chow{\operatorname{A}}

\newcommand\underaut{\operatorname{\underline {Aut}}}

\newcommand\autom[1]{\underaut \bigl(\mathbb{P}^n_{\mathbb{Z}},
\mathcal{O}(#1)\bigr)}

\newcommand\automo[1]{\underaut \bigl(\mathbb{P}^1_{\mathbb{Z}},
\mathcal{O}(#1)\bigr)}

\newcommand\biautomo[2]{\underaut \bigl(\mathbb{P}^1_{\mathbb{Z}},
\mathcal{O}(#1), \mathcal{O}(#2)\bigr)}

\begin{document}

\title[Stacks of cyclic covers of projective spaces]
{Stacks of cyclic covers of projective spaces}

\author{Alessandro Arsie}

\address{Dipartimento di Matematica\\
Universit\`a di Bologna\\
40126 Bologna\\ Italy}
\email{arsie@dm.unibo.it}

\author{Angelo Vistoli}
\email{vistoli@dm.unibo.it}


\begin{abstract}
We define stacks of uniform cyclic covers of Brauer-Severi
schemes, proving that they can be realized as quotient stacks of
open subsets of representations, and compute the
Picard group for the open substacks parametrizing  smooth uniform
cyclic covers. Moreover, we give an analogous description for
stacks parametrizing triple cyclic covers of Brauer-Severi
schemes of  rank $1$ that are not necessarily uniform, and give
a presentation of the  Picard group of the substacks corresponding
to smooth triple cyclic covers.

\end{abstract}

\date{April 7, 2003}

\thanks{Both authors partially supported by the University of
Bologna, funds for selected research topics.}

\maketitle



\section{Introduction}
In \cite{M2} the second author described the stack of
$\mathcal{M}_2$ of smooth curves of genus $2$ as the quotient
stack of an open subscheme of a  representation of $\GL_2$, and
used this description to compute its integral  Chow ring. In
particular, he reproved the known result that its Picard group is
cyclic of order $10$. The key point for the existence of such a
presentation for $\mathcal{M}_2$ is the fact that any smooth
curve of genus $2$ is hyperelliptic.

In this work we define a much wider class of stacks,
parametrizing families of uniform cyclic cover of projective
spaces, that can be
realized as quotient stacks of an open subset of a
representation. Special cases are the stack $\mathcal{M}_2$, the
stacks
$\mathcal{H}_g$ parametrizing hyperelliptic curves of genus $g$
and also the stack parametrizing K3 surfaces
expressed as double covers of $\mathbb{P}^2$ ramified
along a smooth sextic, (up to an automorphism of
$\mathbb{P}^2$). Again, the key idea is that for the objects
involved in families of uniform cyclic covers, one has a concrete
description in terms of polynomials and equations, so that the
corresponding  stack is obtained as a quotient stack of an affine
space parametrizing the corresponding polynomials, modulo the
action of the relevant group.

The paper is organized as follows. In Section~\ref{sec:uniform} we
give the main definitions and constructions for uniform cyclic
covers of a scheme (these are essentially what were known as
simple cyclic covers, see
\cite{C}). A detailed analysis of these and other types of covers
can be found in \cite{Cov}. Moreover we set up the general
categorical framework for uniform cyclic covers over a fixed
scheme.

In Section~\ref{sec:uniform-projective} we restrict our analysis
to uniform cyclic covers of families of projective spaces, i.e.
Brauer-Severi schemes. We introduce our main object of interest,
the fibered categories $\cov$ that parametrize families of
\tots over Brauer-Severi schemes.

In Section~\ref{sec:quotient} we describe $\cov$ and $\covsm$ (the
open substack corresponding to smooth \tots) as quotient
stacks. We also suggest a natural compactification of $\covsm$ via
Kirwan's procedure in Section~\ref{sec:quotient}.

Section~\ref{sec:picard} is dedicated to the computation of the
integral Picard group of the stack $\covsm$; we show that is it
cyclic of order $r(rd-1)^n\mathrm{gcd}(d,n+1)$. As a corollary we
get immediately the the Picard group of the stack
$\mathcal{H}_g$ of hyperelliptic curves of genus $g$ is cyclic of
order $2(2g+1)$ if $g$ is even, and $4(2g+1)$ if $g$ is odd.

Finally, in Section~\ref{sec:triple} we define and study the
stacks $\triple$ of cyclic triple (not necessarily uniform) covers
of the projective line, and its open substack $\triplesm$
corresponding to smooth covers. We prove that this stack can also
be represented as a quotient stack, and we give a presentation of
its Picard group.

\subsection{Acknowledgments} The authors are grateful to Zinovy
Reichstein and Rita Pardini for useful discussions.

This paper can be considered as tribute to the incredible mathematical
insight of David Mumford, who wrote \cite{mum} almost forty years ago.
Reading this paper was a high point in the mathematical education of the
second author.

\section{\Tots of a scheme}\label{sec:uniform}

Fix a positive integer $r$; we will denote by $\mmu = \mmu[r,
\mathbb{Z}]$ the group scheme of $r\th$ roots of $1$ over $\spec
\mathbb{Z}$.

\begin{definition}
Let $Y$ be a scheme. A \emph{\tot of degree $r$} of $Y$ consists
of a morphism of schemes $f \colon X \to Y$ together with an
action of the group scheme $\mmu$ on $X$, such that for each
point $q$ of $Y$, there is an affine neighborhood $V =  \spec R$
of $q$ in $Y$, together with an element $h \in R$ that is not a
zero divisor, and an isomorphism of $V$-schemes $f^{-1}(V) \simeq
\spec R[x]/(x^r - h)$ which is $\mmu$-equivariant, when the right
hand side is given the obvious actions.

\end{definition}

These coverings should be properly called \emph{dual cyclic},
rather than cyclic, as
$\mmu$ is Cartier dual to the constant group scheme $\mathbb{Z}/r
\mathbb{Z}$; but we avoid this, not to make the terminology
unduly heavy. In literature, they are also known as \emph{simple
cyclic cover}.

If $X \to Y$ is a \tot of degree
$r$, then $Y = X/\mmu$, so in fact $Y$ is determined by the action of
$\mmu$ on $X$.

\Tots of a scheme $Y$ form a category, that we
denote by $\covof{Y}$. The arrows are $\mmu$-equivariant
isomorphisms of schemes over $Y$; all the arrows are invertible,
so this category is a groupoid.

There is a very well known description of \tots, as follows. If
$f \colon X
\to Y$ is a
\tot, the sheaf of
$\mathcal{O}_Y$-algebras
$f_* \mathcal{O}_X$ admits an action of $\mmu$, hence there is a
direct sum decomposition
     \[
     f_* \mathcal{O}_X = \mathcal{L}_0 \oplus \mathcal{L}_1 \oplus
    \dots \oplus\mathcal{L}_{r-1},
     \]
where $\mathcal{L}_i$ is the subsheaf of $f_* \mathcal{O}_X$ of
sections
$s$ where the action of $\mmu$ is described by the rule $(t,s)
\mapsto t^i s$. The multiplication is $\mmu$ equivariant,
therefore for each $i = 0$,
\dots,~$r-1$ there is an induced homomorphism
$\mathcal{L}_1^{\otimes i}
\to \mathcal{L}_i$, and also $\mathcal{L}_1^{\otimes r} \to
\mathcal{L}_0$. The local description of the morphism $X
\to Y$ shows that the following facts are true.

\begin{enumeratea}

\item Each $\mathcal{L}_i$ is an invertible sheaf on $Y$.

\item $\mathcal{L}_0 = \mathcal{O}_Y$.

\item For each $i = 0$, \dots,~$r-1$, the homomorphism
$\mathcal{L}_1^{\otimes i} \to \mathcal{L}_i$  is an isomorphism.

\item The homomorphism $\mathcal{L}_1^{\otimes r} \to
\mathcal{O}_Y$ is injective.

\end{enumeratea}

The image of $\mathcal{L}_1^{\otimes r}$ in $\mathcal{O}_Y$ is
the sheaf of ideals of a Cartier divisor on $Y$, that we denote
by $\Delta_f$ or
$\Delta_{X/Y}$, and we call the \emph{branch divisor} of
the \tot. If $V = \spec R$ is an open affine subset of $Y$, such
that $f^{-1}(V)
\simeq \spec R[x]/(x^r - h)$, as in the definition, then the
restriction of $\Delta_f$ to $V$ is the divisor of $h$.

Conversely, assume that we are given a scheme $Y$ with an
invertible sheaf
$\mathcal{L}$, together with an injective homomorphism
$\phi \colon \mathcal{L}^{\otimes r} \to \mathcal{O}_Y$. We can
give the sheaf of $\mathcal{O}_Y$-modules
    \[
    \mathcal{O}_Y \oplus \mathcal{L} \oplus \mathcal{L}^{\otimes 2}
    \oplus \dots \oplus\mathcal{L}^{\otimes (r-1)}
    \]
a structure of $(\mathbb{Z}/r \mathbb{Z})$-graded algebra, by
defining the product of an element $s\in
\mathcal{L}^{\otimes i}$ and $t
\in \mathcal{L}^{\otimes j}$ as
    \[
    s \otimes t \in\mathcal{L}^{\otimes(i + j)}
    \]
if $i+j < r$, and as
    \[
    \phi \otimes \id(s \otimes t) \in
    \mathcal{L}^{\otimes(i + j - r)}
    \]
if $i+j \ge r$,
    \[
    \phi \otimes \id
    \colon  \mathcal{L}^{\otimes(i + j)} \to \mathcal{O}_Y \otimes
    \mathcal{L}^{\otimes(i + j - r)} =
    \mathcal{L}^{\otimes(i + j - r)}
    \]
is the obvious homomorphism. Consider the relative spectrum
$X$ of this sheaf of algebras: the $\mathbb{Z}/r \mathbb{Z}$
grading yields an action of $\mmu$ over $X$, and it is immediate
to verify that in fact $X \to Y$ is a \tot.

This analysis leads to the following conclusion. Define a category
$\altcovof Y$, whose objects $(\mathcal{L}, \phi)$ are invertible
sheaves
$\mathcal{L}$ on $Y$, together with an injective homomorphism of
$\mathcal{O}(Y)$-modules $\phi\colon
\mathcal{L}^{\otimes r} \to \mathcal{O}_Y$. The arrows $\alpha
\colon (\mathcal{L}, \phi) \to (\mathcal{M}, \psi)$ are
isomorphisms of invertible sheaves $\alpha \colon \mathcal{L}
\simeq \mathcal{M}$, making the diagram
     \[
     \xymatrix@-10pt{
     \mathcal{L}^{\otimes r} \ar[rr]^{\alpha ^{\otimes r}}
    \ar[rd]_\phi
     && \mathcal{M}^{\otimes r}\ar[ld]^\psi\\
     & \mathcal{O}_Y
     }
     \]
commutative.

\begin{proposition}
There is an equivalence of categories between the category $\covof Y$ and
the category $\altcovof Y$.
\end{proposition}

Given a \tot $f \colon X \to Y$, the pullback of $\Delta_f$ to
$X$ is a Cartier divisor, which is of the form $r D_f$, where
$D_f$ is a Cartier divisor on $X$, whose sheaf of ideals is the
pullback $f^* \mathcal{L}$, where $\mathcal{L}$ is the invertible
sheaf associated with $f \colon X
\to Y$. The restriction $D_f \to \Delta_f$ is an isomorphism.

There is a problem with defining pullbacks of \tots: if $f \colon
X \to Y$ is a \tot, and $Y' \to Y$ a morphism of schemes, the
pullback $X' \eqdef Y'
\times_Y X$ acquires a natural actions of $\mmu$, but the
projection $f'
\colon X' \to Y'$ is a \tot if and only if the pullback of the
branch divisor $\Delta_f$ to $Y'$ is still a Cartier
divisor. This problem does not arise in a relative context, which
the one that we are interested in.

\begin{definition}
Let $Y \to S$ be a morphism of schemes. A \emph{relative \tot} $f
\colon X
\to Y$ is a \tot, such that the branch divisor $\Delta_f$
is flat over $S$.
\end{definition}

By the local criterion of flatness, $f \colon X \to Y$ is a
relative \tot if and only if $\Delta_f$ remains a Cartier divisor
when restricted to any of the fiber of $Y \to S$.

The relative \tots over $Y \to S$ form a full subcategory of
$\covof Y$, denoted by $\covof{Y/S}$

If $f \colon X \to Y$ is a relative \tot over $Y \to S$, and $S'
\to S$ is an arbitrary morphism of schemes, then the pullback of
$\Delta_f$ to $S'
\times_S Y$ is still a Cartier divisor, so the projection $S'
\times_S X
\to S' \times_S Y$ is a relative \tot.

\begin{definition}\label{def:smooth}

A relative \tot $f \colon X \to Y$ over a morphism $Y \to S$ is
\emph{smooth over $S$} if both $Y$ and the branch divisor
$\Delta_f$ are smooth over $S$.
\end{definition}

The proof of the following is straighforward.

\begin{proposition}\label{prop:describe-smooth}
Let $Y \to S$ be a smooth morphism, $f \colon X \to Y$ a relative \tot of
degree $r$. Then $f$ is
a smooth \tot over $S$ if and only if $X$ is smooth over $S$.
\end{proposition}

\section{\Tots of projective spaces}\label{sec:uniform-projective}

We are interested in relative \tots $f \colon X \to P$ of degree
$r$, where
$P \to S$ is a Brauer--Severi scheme. Given such a thing,
consider the invertible sheaf
$\mathcal{L}$ of sections of $f_* \mathcal{O}_X$ on which $\mmu$
acts via multiplication. The degree of such a invertible sheaf on the
geometric fibers of $P \to S$ is a local invariant. We say that
such a \tot
\emph{has branch degree $d$} if the degree of $\mathcal{L}$
is $d$ on every fiber; the degree of the branch divisor is
then equal to
$rd$ (so perhaps this is not great terminology).

Fix three positive integers $n$, $r$ and $d$. We are interested
in the category $\cov$, defined as follows.

An object $(X \stackrel f\to P \to S)$ of $\cov$ is a relative
\tot $f
\colon X \to P$ of degree $r$ and branch degree $d$, where
$P \to S$ is a Brauer--Severi scheme of relative dimension $n$.

An arrow from $(X' \stackrel{f'}\to P' \to S')$ to $(X \stackrel
f\to P
\to S)$ is a commutative diagram
    \[
    \xymatrix{
    X' \ar[r]^{f'}\ar[d] &P'\ar[r]\ar[d] &S'\ar[d] \\
    X \ar[r]^{f}   &P\ar[r]  &S
    }
    \]
where both squares are cartesian, and the left hand column is
$\mmu$-equivariant.

We can reformulate the definition as follows.

\begin{proposition}\label{prop:alt-description}
The category $\cov$ is equivalent to the category $\altaltcov$
defined as follows. The objects are flat and proper morphisms $X
\to S$ of schemes, together with an action of $\mmu$ on $X$
leaving $X \to S$ invariant, satisfying the following condition:
for any geometric point $s \colon
\spec \Omega \to X$, the action on $\mmu$ on the geometric fiber
$X_s$ is faithful, the quotient
$X_s/\mmu$ is isomorphic to $\mathbb{P}^n_{\spec\Omega}$, and the
projection $X_s \to X_s/\mmu$ makes $X_s$ into a \tot of
$X_s/\mmu$, with degree $r$ and branch index $d$.

The arrows from $X' \to S'$ to $X \to S$ are commutative squares
     \[
     \xymatrix{
     X'\ar[d]\ar[r] & X\ar[d]\\
     S'\ar[r] & S
     }
     \]
such that the top row is $\mmu$-equivariant.
\end{proposition}

\begin{proof}

Given an object $(X \to P \to S)$ of $\cov$, we have that the
composition
$X \to S$ gives an object of $\altaltcov$; this, together with the
analogous construction for arrows, defines a functor $\cov \to
\altaltcov$. To go in the other direction we need a lemma.

\begin{lemma}\label{lem:basechange}
If $X \to S$ is a morphism of schemes, and there is given an
action of
$\mmu$ on $X$ leaving $X \to S$ invariant, then the formation of
the quotient $X/\mmu$ commutes with base change on $S$.
Furthermore, if $X$ is flat over $S$, so is $X/\mmu$.
\end{lemma}

\begin{proof}
Both part of the statement are standard consequences of the
fact that
$\mmu$ is a diagonalizable group scheme over $\spec \mathbb{Z}$.
\end{proof}

Suppose that $X \to S$ is an object of $\altaltcov$, and factor it
as $X \to P \to S$, where $P = X/\mmu$. Obviously $P$ is proper
over $S$. The lemma implies that it is also flat over $S$, and
that the geometric fibers are projective spaces; hence, by a well
known theorem of Grothendieck, $P$ is a Brauer--Severi scheme
over $S$.

Also, the restrictions of the projection morphism $f \colon X \to P$ over the
points of $S$ is flat, so, by the local criterion of flatness, $f$ itself is
flat. It is also finite, so
$X$ can be thought of as the relative spectrum on $P$ of the
locally free sheaf of algebras $f_* \mathcal{O}_X$, which we can
decompose as
     \[
     f_* \mathcal{O}_X = \mathcal{L}_0 \oplus \mathcal{L}_1 \oplus \dots
     \oplus\mathcal{L}_{r-1},
     \]
using the action of $\mmu$. For each $i = 0$, \dots,~$r-1$, the
natural homomorphism
$\mathcal{L}_1^{\otimes i} \to \mathcal{L}_i$ is an isomorphism
on each geometric fiber, hence it is an isomorphism; furthermore
$\mathcal{L}_1^{\otimes r} \to \mathcal{O}_P$ is injective on the
geometric fibers. This means that $f \colon X \to P$ is a \tot,
hence $X \to P \to S$ is an object of $\cov$.

It is very easily checked that this extends naturally to a functor
$\altaltcov \to \cov$, and this gives a quasi-inverse to the functor
above.
\end{proof}

\begin{remark}\label{rmk:altcov}
It is also convenient to define a fibered category $\altcov$,
in which an object over a scheme $S$ consists of the following
set of data: a Brauer-Severi  scheme $P \to S$, an invertible
sheaf
$\mathcal{L}$ on $P$, which restricts to a invertible sheaf of
degree
$-d$ on any geometric fiber, and an injection
$i\colon \mathcal{L}^{\otimes r} \to \mathcal{O}_P$, which
remains injective when restricted to any geometric fiber. The
morphisms are defined in the  obvious way. Clearly there is a
morphism of fibered categories $p\colon
\cov \to \altcov$ sending the object $(X \to P \to S)$ to the
triple  ($P
\to S$, $\mathcal{L}$, $i\colon \mathcal{L}^{\otimes r} \to
\mathcal{O}_P$) and acting in the obvious way on morphisms.  This
correspondence is also
   an equivalence of fibered category, as it is immediate to see,
since $X$ can  be recovered as
$\underline{\spec}_{\mathcal{O}_P}(\mathcal{O}_P\oplus
\mathcal{L}\oplus \dots \oplus \mathcal{L}^{\otimes r-1})$.

\end{remark}

We denote by $\covsm$ the full subcategory of $\cov$ consisting of
relative \tots $X \to P \to S$ which are smooth over the base.

There is a natural forgetful functor from $\cov$ to the category
of schemes, sending $(X \stackrel f\to P \to S)$ to $S$; this
makes $\cov$ into a fiber category over the category of schemes,
and $\covsm$ is a fibered subcategory.

  From now on, if $R$ is a commutative ring, we will write
$\cov_R$ for the fiber product of $\cov$ with the category of
schemes over $R$; the objects of $\cov_R$ are pairs $\bigl((X \to
P \to S), S \to \spec R\bigr)$ consisting of an object of $\cov$
and of a morphism of schemes. The arrows are defined in the
obvious way. There will be obvious variant of this notation, such as
$\covsm_R$ and $\altcov_R$ (the category $\altcov$ is defined
above).

The category
$\covsm_{\mathbb{Z}[1/r]}$ is fibered on the category of schemes over
$\spec\mathbb{Z}[1/r]$. It has a simple description, using the
equivalent description of $\altcov$ given in
Proposition~\ref{prop:alt-description}.

\begin{proposition}
The fibered category $\covsm_{\mathbb{Z}[1/r]}$ is equivalent to
the full subcategory of $\altcov$ consisting of objects $X \to S$
which are smooth as morphisms of schemes, and where $S$ is a
scheme over
$\spec\mathbb{Z}[1/r]$.
\end{proposition}

\begin{proof}
This follows from Proposition~\ref{prop:alt-description} and
Proposition~\ref{prop:describe-smooth}.
\end{proof}

Here are some examples of our construction.

\begin{examples}\hfil\call{ex:examples}
\begin{enumeratea}

\item For each $g \ge 2$, the fibered category
$\covsm[1,2,g+1]_{\mathbb{Z}[1/2]}$ is a closed substack of the
stack
$\mathcal{M}_g$ of smooth curves of genus $g$, whose geometric
points are the hyperelliptic curves.

In particular, $\covsm[1,2,3]_{\mathbb{Z}[1/2]}$ coincides with
$\mathcal{M}_2$.

\itemref{1}  We do not know if the category $\covsm[1,2,2]$ has
appeared in the literature before. Its objects are smooth
families $X \to S$ of curves of genus $1$ over a scheme on
$\spec\mathbb{Z}[1/2]$, together with an effective divisor
$\Sigma \subseteq X$, such that the restriction $\Sigma \to S$ is
\'etale of degree 4, and that $\Sigma$ is invariant under the
action of the
$2$-torsion part $\sideset{_2}{^0}{\pic}(X/S) \to S$ of the
associated elliptic curve.

\item Consider the category $\covsm[2,2,3]$ of double covers of a
projective plane, ramified over a smooth sextic curve. In
characteristic different from $2$, the resulting surfaces are K3
surfaces, of a special and well studied type.

\end{enumeratea}
\end{examples}

\begin{remark}\label{rmk:why-not1}
More generally, one might be interested in flat morphisms $f
\colon X \to P$, where $P \to S$ is a Brauer--Severi scheme, together
with an action of
$\mmu$ on $X$ leaving $f$ invariant, such that there exists an
open subscheme $U$ of $P$, dense in every fiber of $P \to S$, such that, over
$U$, the restriction of $f$ is a $\mmu$-torsor. Among these,
uniform cyclic covers are special in two ways.

First of all, they are totally ramified (that is, the action of
$\mmu$ is free outside of the fixed locus); of course, this is
only a restriction when $r$ is not a prime.

Also, the action of $\mmu$ around a fixed point is of very
restricted type; for example, if we are looking at a smooth \tot
$f
\colon X \to
\mathbb{P}^n$ defined over $\mathbb{C}$, then the restriction $X'
\to
\mathbb{P}^n \setminus \Delta_f$ is a Galois covering with group
$\mmu$, whose restriction to a small loop $L \simeq \mathbb{S}^1$
around a smooth point of $\Delta_f$ corresponds to the canonical
generator of $\H^1(L,
\mmu) = \mathbb{Z}/r \mathbb{Z}$.

If $n>1$, one might consider this not to be a serious restriction;
for example, if $r$ is a prime power, $f \colon X \to
\mathbb{P}^n$ is a flat morphism defined over a field of
characteristic prime to $r$, and there is an action of $\mmu$
on $X$ leaving $f$ invariant, such that generically $X$ is a
torsor over $\mathbb{P}^n$, and
$X$ is smooth over the base field, then it is not hard to show
that we can make $f \colon X \to \mathbb{P}^n$ into a \tot by
changing the action by an automorphism of $\mmu$; thus the
resulting stack is a disjoint union of copies of $\covsm$.

When $r$ is not a prime power, then this is not true anymore;
however, one can still describe this stack as an open substack of
products of stacks of type $\covsm[n,r_i,d_i]$.

Things are altogether different when $n = 1$ and $r>2$; here the
branch divisor will almost never be irreducible, and cyclic
coverings of $\mathbb{P}^1$ that are not uniform are very common.
We will describe the situation for $\mmu[3]$-covers in
Section~\ref{sec:triple}.
\end{remark}

\begin{remark} The stack $\cov$ itself is not particularly
useful; the objects involved are highly unstable. We will be
mostly interested in $\covsm$; there is a natural
compactification of it, via Kirwan's procedure, as explained in
Remark~\ref{rmk:Kirwan}.
\end{remark}

\section{$\cov$ as a quotient stack}\label{sec:quotient}

For each triple $n$, $r$ and $d$, consider the space $\forms$ of
homogenous forms of degree $rd$ in $n+1$ indeterminates; we can
think about
$\forms$ as the spectrum of the polynomial ring
$\mathbb{Z}[a_I]$, where
$a_I$ is an indeterminate, and $I$ varies over the set of
functions $I
\colon
\{0, \dots, n\} \to \mathbb{N}$ with $\sum_k I(k) = rd$, so
$\forms$ is an affine space of dimension $\binom{rd+n}{n}$ over
$\mathbb{Z}$.

We also write $\pforms$ for the projective space of lines in
$\forms$ (in this context this convention seems more natural than
Grothendieck's).

We denote by $\formsnoto$ the complement of the zero section
$\spec\mathbb{Z} \into
\forms$, and by $\formssm \subseteq \formsnoto$ the open subscheme
corresponding to smooth forms.

There is a natural action of $\GL_{n+1} = \GL_{n+1, \mathbb{Z}}$
on
$\forms$, defined, in functorial notation, by $A \cdot f(x) =
f(A^{-1}x)$. The subgroup scheme $\mmu[d] \subseteq \GL_{n+1}$,
embedded by sending a $d\th$ root of $1$ $\alpha$ into the
diagonal matrix $\alpha I_{n+1}$, acts trivially on $\forms$, so
this induces an action of the quotient
$\GL_{n+1}/\mmu[d]$ on $\forms$, leaving the open subschemes
$\formsnoto$ and $\formssm$ invariant.

\begin{theorem}\label{thm:quotient}
The fibered category $\cov$ is isomorphic to the quotient stack
     \[
     [\formsnoto/(\GL_{n+1}/\mmu[d])]
     \]
by the action described above.

Furthemore, if $R$ is a commutative ring, and $F \in
\formsnoto(R)$ a form of degree~$rd$ whose coefficients generate
the trivial ideal $R$, the branch divisor $\Delta_f
\subseteq \mathbb{P}^n_R$ of the associated \tot
$f \colon X \to \mathbb{P}^n_R$ is the hypersurface of
$\mathbb{P}^n_R$ defined by $F$.
\end{theorem}

\begin{proof}

To prove the theorem, we identify $\cov$
with $\altcov$, the fibered category of Remark~\ref{rmk:altcov}.

Consider the auxiliary fibered category $\covv$, whose objects
over a base scheme $S$ are given as pairs consisting of an object
$\bigl(P \to S, \mathcal{L}, i \colon \mathcal{L}^{\otimes r} \to
\mathcal{O}_P\bigr)$ in $\cov(S)$, plus an isomorphism
$\phi\colon (P,\mathcal{L}) \simeq (\mathbb{P}^n_S,
\mathcal{O}(-d))$, over $S$ (by this we mean the pair consisting
of an isomorphism of $S$-schemes
$\phi_0  \colon P \simeq \mathbb{P}^n_S$, plus an isomorphism
$\phi_1 \colon \mathcal{L} \simeq \phi_0^*\mathcal{O}(-d)$). The
arrows in
$\covv$ are arrows in $\cov$ preserving the isomorphisms $\phi$.

The obvious projection
from $\covv$ to the category of schemes makes it into a
category fibered in groupoids. In fact, no object of $\covv$ has
a nontrivial automorphism mapping to identity in the category of
schemes, so $\covv$ is equivalent to a functor. One has a morphism
of fibered categories from $\covv$ to $\altcov$ by
forgetting the isomorphism $\phi$.

Let us define a base-preserving functor from $\covv$ to
$\formsnoto$. For any object of
$\covv(S)$ take the composition
    \[
    \phi \circ i \circ(\phi^{-1})^{\otimes r}\colon
    \mathcal{O}_{\mathbb{P}^n_S}(-rd) \to
    \mathcal{O}_{\mathbb{P}^n_S},
    \]
corresponding to a section of
$\mathcal{O}_{\mathbb{P}^n_S}(rd)$ that does not vanish on any
fiber of $\mathbb{P}^n_S \to S$, that is, to an element of
$\formsnoto(S)$. There is also a base preserving functor in the
other direction, by sending a section $f \in
\mathcal{O}_{\mathbb{P}^n_S}(rd)$, thought of as a homomorphism $f
\colon \mathcal{O}_{\mathbb{P}^{n}_S}(-rd) \to
\mathcal{O}_{\mathbb{P}^n_S}$, into the object
    \[
    \bigl( \mathbb{P}^n_S \to S, \mathcal{O}(-d), f \colon
    \mathcal{O}(-d)^{\otimes r} \to \mathcal{O}, \id \colon
    (\mathbb{P}^n_S, \mathcal{O}(-d)) \to
    (\mathbb{P}^n_S, \mathcal{O}(-d))
    \bigr)
    \]
of $\covv (S)$. It is straightfoward to check that this gives a
quasi-inverse to the previous functor; so we get an equivalence
of $\covv$ with $\formsnoto$.

Now, for each interger $e$ consider the functor $\autom{e}$ from
schemes to groups sending each scheme $S$ into group of
automorphisms of the pair
$\bigl(\mathbb{P}^n_S, \mathcal{O}(e)\bigr)$ over the identity
on $S$. This is a sheaf in the fppf topology. Clearly,
$\autom{1}$ can be identified with $\GL_{n+1, \mathbb{Z}}$; an
isomorphism of the pair
$(\mathbb{P}^n_S,\mathcal{O}(1))$, gives via $\pi\colon
\mathbb{P}^n_S \to S $ an automorphism of
$\pi_{*}\mathcal{O}(1)=\mathcal{O}^{n+1}_S$ as an
$\mathcal{O}_S$-module, and conversely. There is a natural
homomorphism of sheaves of groups
    \[
   \autom{1} \to \autom{e}
    \]
sending each automorphism $(\phi_0, \phi_1) \colon
\bigl(\mathbb{P}^n_S, \mathcal{O}(1)\bigr) \simeq
\bigl(\mathbb{P}^n_S, \mathcal{O}(1)\bigr)$ into
    \[
    (\phi_0, \phi_1^{\otimes e}) \colon
    \bigl(\mathbb{P}^n_S, \mathcal{O}(e)\bigr) \simeq
    \bigl(\mathbb{P}^n_S, \mathcal{O}(1)\bigr).
    \]
It is easy to check that this is a surjective homomorphism of
fppf sheaves. If we identify $\autom{1}$ with $\GL_{n+1,
\mathbb{Z}}$, then the kernel of this homomorphism is the
subgroup $\mmu[|e|, \mathbb{Z}]$ embedded diagonally. So we get
an isomorphism
    \[
    \autom{-d} \simeq
    \GL_{n+1, \mathbb{Z}}/\mmu[d, \mathbb{Z}]
    \]

There is a left action of $\autom{-d}$ on $\covv$; if
    \[
    \bigl(P \to S, \mathcal{L}, i \colon \mathcal{L}^{\otimes r}
    \to \mathcal{O}_P, \phi\colon (P,\mathcal{L}) \simeq
    (\mathbb{P}^n_S, \mathcal{O}(-d))\bigr)
    \]
is an object of $\covv(S)$, and
    \[
    \alpha \colon \bigl(\mathbb{P}^n_S, \mathcal{O}(-d)\bigr)
    \simeq \bigl(\mathbb{P}^n_S, \mathcal{O}(-d)\bigr)
    \]
an element of
$\autom{-d}$, we associate with these the object
    \[
    \bigl(P \to S, \mathcal{L}, i \colon \mathcal{L}^{\otimes r}
    \to \mathcal{O}_P, \alpha  \circ  \phi \colon (P,\mathcal{L})
    \simeq(\mathbb{P}^n_S, \mathcal{O}(-d))\bigr).
    \]

Furthermore, given an invertible sheaf $\mathcal{L}$ on
$P \to S$ whose degree is $-d$ on every geometric fiber, there is an
fppf covering $S' \to S$, such that the pullback of the pair $(P,
\mathcal{L})$ to $S'$ is isomorphic to
$(\mathbb{P}^n_{S'},\mathcal{O}(-d))$; this
fact, plus descent theory, implies that the forgetful morphism
$\covv \to \cov$ makes $\covv$ into a principal bundle with group
$\autom{-d} = \GL_{n+1, \mathbb{Z}}/\mmu[r, \mathbb{Z}]$.

If we identify $\covv$ with $\formsnoto$, we obtain that $\cov$ is
isomorphic to the quotient stack $[\formsnoto/(\GL_{n+1,
\mathbb{Z}}/\mmu[r, \mathbb{Z}])]$. We only have left to identify
the action explicitly. But from the description above it is easy
to check that $\GL_{n+1} = \autom{1}$ acts by the usual action
$(f \cdot A) (x) = f(A^{-1}x)$, so the action of its quotient is
the one described above.

The last statement follows easily by construction.
\end{proof}

The following corollary is a direct application of
Theorem~\ref{thm:quotient}.

\begin{corollary}\label{cor:describe-covsm}
The fibered category $\covsm$ is equivalent to the quotient stack
     \[
     [\formssm/(\GL_{n+1}/\mmu[d])]
     \]
by the action of described above.
\end{corollary}

In particular, $\cov$ is an irreducible smooth algebraic stack of finite
type over
$\spec\mathbb{Z}$, of relative dimension
     \[
     \binom{rd + n}{n} - (n+1)^2,
     \]
and $\covsm$ is an open substack, hence it also smooth of the same
dimension So, for example, the dimension of the stack of
hyperelliptic curves
$\covsm[1,2,g+1]$ is $2g-1$, as it should, and the dimension of
the stack of K3 surfaces $\covsm[2,2,3]$ is 19.

The fact that hypersurfaces of degree at least three are stable
for the action of $\SL_{n+1}$ implies that when $d > 1$ the
diagonal of $\covsm$ is finite, and its moduli space is
quasiprojective over $\spec\mathbb{Z}$. Also, again for
$d > 1$ the restriction of
$\covsm_{\spec\mathbb{Z}[1/rd]}$ is a Deligne--Mumford stack over
$\spec\mathbb{Z}[1/rd]$.

\begin{remark}\label{rmk:Kirwan}
Assume that $d$ is at least $3$. Then, if we look at the natural
action of
$\SL_{n+1}$ on the projectivization $\pforms$ of $\forms$, the
points corresponding to smooth hypersurfaces are stable. This
implies that we can apply Kirwan's procedure (see \cite{Kir}), to
get a canonical
$\GL_{n+1}/\mmu[d]$-equivariant morphism $K(n,rd) \to
\formsnoto$, which is an isomorphism over $\formssm$, such that
the action of $\GL_{n+1}/\mmu[d]$ is proper, and the geometric
quotient $K(n,rd)/(\GL_{n+1}/\mmu[d])$ is a projective scheme
over $\spec\mathbb{Z}$. The quotient stack
    \[
    \covcomp = [K(n,rd)/(\GL_{n+1}/\mmu[d])]
    \]
is an Artin stack with
finite diagonal and projective moduli space, yielding a canonical
compactification of $\covsm$; this seems like a much more natural
object than $\cov$.

One could try to investigate the stacks $\covcomp$, and in
particular describe their objects directly. This seems very
complicated in dimension higher than 2, but at least for
$n=1$ the problem should be approachable. If we exclude
characteristic~$2$ then $\cov[1,2,3]$ is the stack
$\mathcal{M}_2$ of smooth curves of genus~2, and one can check
that
$\covcomp[1,2,3]$ is not isomorphic to the stack
$\overline{\mathcal{M}}_2$ of stable curves of genus~2, although
it would seem that it gives the same moduli space. But in the next
case, $\cov[1,2,4]$ is the stack of smooth hyperelliptic curves of
genus~3, and one can easily see that
$\covcomp[1,2,4]$ does not coincide with the closure of
$\cov[1,2,4]$ inside
$\overline{\mathcal{M}}_3$, not even at the level of moduli
spaces; so the stack of hyperelliptic curves of fixed genus $g$
has two natural compactifications, and in general they do not
coincide.

It would be interesting to investigate these two
compactifications, and try to determine if they have any
relations.
\end{remark}

The group $\GL_{n+1}/\mmu[d]$ appearing in the statement of the
theorem can sometimes be written in a more familiar form. The
following is straightforward.

\begin{proposition} \label{prop:isom-groups}\hfil
\begin{enumeratea}

\item If $d \equiv 0 \pmod{n+1}$, write $d = q(n+1)$. The
homomorphism of group schemes over $\mathbb{Z}$
     \[
     \GL_{n+1}/\mmu[d] \longrightarrow \gm \times \PGL_{n+1}
     \]
defined by
     \[
     [A] \mapsto (\det(A)^q, [A])
     \]
is an isomorphism.

\item If $d \equiv 1 \pmod{n+1}$, write $d = q(n+1) + 1$. The
homomorphism of group schemes over $\mathbb{Z}$
     \[
     \GL_{n+1}/\mmu[d] \to \GL_{n+1}
     \]
defined by
     \[
     [A] \mapsto \det(A)^q A
     \]
is an isomorphism.

\item If $d \equiv -1 \pmod{n+1}$, write $d = q(n+1) - 1$. The
homomorphism of group schemes over $\mathbb{Z}$
     \[
     \GL_{n+1}/\mmu[d] \to \GL_{n+1}
     \]
defined by
     \[
     [A] \mapsto \det(A)^{-q} A
     \]
is an isomorphism.

\end{enumeratea}
\end{proposition}

\begin{remark}
One can show that the group scheme $\GL_{n+1}/\mmu[d]$ is
isomorphic to
$\GL_{n+1}$ if and only if $d \equiv \pm 1 \pmod{n+1}$; on the
other hand
$\GL_{n+1}/\mmu[d]$ is special (in the sense that every
$\GL_{n+1}/\mmu[d]$-torsor is locally trivial in the Zariski
topology) if and only if $d$ is prime to
$n+1$ (Zinovy Reichstein pointed this out to us). Experience
teaches that special groups are infinitely easier to handle than
nonspecial ones; so, computing basic invariants of the spaces
$\cov$ and
$\covsm$ (such as Chow rings and cohomology) shoud be much easier
when
$d$ is prime to $n+1$. For this purpose, it would be useful to
gather information about the cohomology and the Chow ring of the
classifying spaces of these groups.
\end{remark}

If we rewrite the action of Theorem~\ref{thm:quotient} via the
isomorphisms of Proposition~\ref{prop:isom-groups} we obtain the
following.

\begin{corollary}\hfil
\begin{enumeratea}

\item If $d \equiv 0 \pmod{n+1}$, write $d = q(n+1)$. Then $\cov$ is
equivalent to the quotient stack
     \[
     [\forms/(\gm \times \PGL_{n+1})]
     \]
by the action defined by the formula
     \[
     (\alpha, [A]) \cdot f(x) = \alpha^{-r}\det(A)^{rq} f(A ^{-1}
    x).
     \]

\item If $d \equiv 1 \pmod{n+1}$, write $d = q(n+1) + 1$. Then
$\cov$ is equivalent to the quotient stack
     \[
     [\forms/(\GL_{n+1})]
     \]
by the action defined by the formula
     \[
     A \cdot f(x) = \det(A)^{rq} f(A ^{-1} x).
     \]

\item If $d \equiv -1 \pmod{n+1}$, write $d = q(n+1) - 1$. Then
$\cov$ is equivalent to the quotient stack
     \[
     [\forms/(\GL_{n+1})]
     \]
by the action defined by the formula
     \[
     A \cdot f(x) = \det(A)^{-rq} f(A ^{-1} x).
     \]

\end{enumeratea}

\end{corollary}

In particular we get the following description of the stack of
hyperelliptic curves.

\begin{corollary}
The stack $\covsm[1,2,g+1]$ of smooth hyperelliptic curves of
genus~$g$ is isomorphic to
\begin{enumeratea}

\item the quotient of $\formssm[1,2g+2]$ by the action of
$\GL_2$ defined by $A \cdot f(x) = \det(A)^{g} f(A ^{-1} x)$ if
$g$ is even, and

\item the quotient of $\formssm[1,2g+2]$ by the action of $\gm
\times\PGL_2$ defined by $(\alpha, [A]) \cdot f(x) =
\alpha^{-2}\det(A)^{g+1} f(A ^{-1} x)$ if $g$ is odd.

\end{enumeratea}
\end{corollary}

When $g = 2$ we recover the description of the stack
$\mathcal{M}_2$ of smooth curves of genus~2 given in \cite{M2};
the derivation here is much simpler, but the method of \cite{M2}
has some independent interest.

\section{Picard groups of stacks of smooth
cyclic coverings}\label{sec:picard}

We will use the description of $\covsm$ given in
Corollary~\ref{cor:describe-covsm} to compute its Picard group,
away from some bad characteristics.

Recall that if $\mathcal{X}$ is an algebraic stack over a scheme $S$, its
Picard group is the group of isomorphism classes of invertible sheaves on
$\mathcal{X}$, with the operation given as usual by tensor product. An
invertible sheaf is a quasicoherent sheaf over $\mathcal{X}$, defined as in
\cite{lm}, which is locally free of rank $1$ when restricted to an atlas.

The Picard group of the stack $\mathcal{M}_{1,1}$ of elliptic curves was
first computed by Mumford in the legendary paper \cite{mum}, written before
the notion of algebraic stack was introduced.

\begin{theorem}\label{thm:picard}
Let $R$ be a unique factorization domain such that the
caracteristic of its quotient field does not divide $2rd$. Then
the Picard group of the stack
$\covsm_{R}$ is cyclic, of order
     \[
     r(rd - 1)^n \gcd(d, n+1).
     \]
\end{theorem}

\begin{proof}

First of all, it follows from the following lemma that we can
assume that $R$ is a field.

\begin{lemma}\label{lem:Picard-to-generic}
Let $\mathcal{X}$ be a flat regular algebraic stack of finite type
over a unique factorization domain $R$ with quotient field $K$.
Assume that the fibers of $\mathcal{X}$ over the closed points of
$\spec R$ are integral. Then the restriction homomorphism
    \[
    \pic \mathcal{X} \longrightarrow \pic (\spec K \times_{\spec
    R} \mathcal{X})
    \]
is an isomorphism.
\end{lemma}

\begin{proof}
The group of divisors $\Div\mathcal{X}$ is the free
abelian group generated by integral closed substacks of
codimension $1$ in $\mathcal{X}$. Effective divisors are defined
in the usual fashion.

The group $\Div\mathcal{X}$ can also be defined as follows:
closed substacks of $\mathcal{X}$ that are local complete
intersection of codimension $1$ form a monoid with the
cancellation property, the operation being defined by taking
products of sheaves of ideals. It is the free abelian monoid on
the set of integral closed substacks of
codimension $1$ in $\mathcal{X}$. The group
$\Div\mathcal{X}$ is the group of quotients of this monoid.

If $f \colon \mathcal{X}' \to \mathcal{X}$ is a dominant morphism
of noetherian regular algebraic stacks, any closed local complete
intersection substack of $\mathcal{X}$ of codimension $1$ pulls
back to a closed local complete
intersection substack of $\mathcal{X}'$ of codimension $1$; this
induces a group homomorphism $f^* \colon \Div \mathcal{X} \to
\Div \mathcal{X}'$.

If $\mathcal{D}$ is a divisor on
$\mathcal{X}$, we can associate with it a divisor $\mathcal{D}_U$
for each smooth morphism $U \to \mathcal{X}$, where $U$ is a
scheme, so that, given two smooth
morphisms $U \to \mathcal{X}$ and $V \to \mathcal{X}$, the
pullbacks of $\mathcal{D}_U$ and $\mathcal{D}_V$ to $U
\times_{\mathcal{X}} V$ coincide with $\mathcal{D}_{U
\times_{\mathcal{X}} V}$. This is done as follows: write
$\mathcal{D}$ as $\mathcal{D}^{+} - \mathcal{D}^{-}$, where
$\mathcal{D}^{+}$ and $\mathcal{D}^{-}$ are effective and do not
intersect in codimension $1$; they correspond to
closed substacks of
$\mathcal{X}$; they pull back to effective divisors
$\mathcal{D}^{+}_U$ and
$\mathcal{D}^{-}_U$ on $U$. We define $\mathcal{D}_U$ to be the
difference $\mathcal{D}^{+}_U - \mathcal{D}^{-}_U$.

Conversely, if we are given a collection of divisors $D_U$ on $U$
for each smooth morphism $U \to
\mathcal{X}$, where $U$ is a scheme, such that given two smooth
morphisms $U \to \mathcal{X}$ and $V \to \mathcal{X}$, the
pullbacks of $D_U$ and $D_V$ to $U \times_{\mathcal{X}} V$
coincide with $D_{U \times_{\mathcal{X}} V}$, there is a unique
divisor $\mathcal{D}$ such that $\mathcal{D}_U = D_U$ for each
smooth morphism $U \to \mathcal{X}$. If we write $D_U = D^+_U -
D^-_U$, where $D^+_U$ and $D^-_U$ are effective and do not
intersect in codimension $1$, for each pair of
smooth morphisms $U \to \mathcal{X}$ and $V \to \mathcal{X}$ we
have $D^+_{U\times_{\mathcal{X}} V} = D^+_U$ and
$D^-_{U\times_{\mathcal{X}} V} = D^-_U$; so the $D^+_U$ and
$D^-_U$ descend to closed substacks of codimension
$1$, $\mathcal{D}^{+}$ and
$\mathcal{D}^{-}$ of $\mathcal{X}$, whose ideals are
locally generated by one element. We set $\mathcal{D} =
\mathcal{D}^{_+} - \mathcal{D}^{-}$.

If $\mathcal{D}$ is a divisor on $\mathcal{X}$, we can associate
with it an invertible sheaf $\mathcal{O}(\mathcal{D})$ on
$\mathcal{X}$, together with a non-vanishing section defined over
the complement of the support of $\mathcal{D}$. Consider the
invertible sheaf
$\mathcal{O}(\mathcal{D}_U)$ defined over $U$ for each smooth
morphism $U
\to \mathcal{X}$. If $U \to \mathcal{X}$ and $V \to \mathcal{X}$
are smooth morphisms, there is a natural isomorphism of
$\mathcal{O}(\mathcal{D}_{U \times_{\mathcal{X}} V})$ with the
pullback of $\mathcal{O}(\mathcal{D}_U)$; these isomorphisms
define the descent data for an invertible sheaf on
$\mathcal{X}$, that we call
$\mathcal{O}(\mathcal{D})$. On the complement of the support of
$\mathcal{D}$ this invertible sheaf is canonically trivial.

This defines a group homomorphism $\Div\mathcal{X} \to \pic
\mathcal{X}$. If  $f \colon \mathcal{X}' \to \mathcal{X}$ is a
dominant morphism of noetherian regular algebraic stacks and
$\mathcal{D}$ is a divisor on $\mathcal{X}$, then
$\mathcal{O}(f^* \mathcal{D})$ is canonically isomorphic to
$f^*\mathcal{O}( \mathcal{D})$.

Conversely, if $\mathcal{L}$ is an invertible sheaf on
$\mathcal{X}$, and $s$ is a nowhere vanishing section of
$\mathcal{L}$ on an open dense substack $\mathcal{U}$, we can
associate with it a divisor $\mathrm{Z}(s)$ on $\mathcal{X}$.
If $\phi \colon U \to \mathcal{X}$ is a smooth morphism, we define
$\mathrm{Z}(s)_U$ to be the divisor of the rational section
$\phi^* s$ of the invertible sheaf $\phi^*\mathcal{L}$ on $U$.

One checks immediately that $s$ extends to a nowhere vanishing
function of the invertible sheaf $\mathcal{L}\otimes
\mathcal{O}(-\mathrm{Z}(s))$; therefore $\mathcal{L}\otimes
\mathcal{O}\bigl(-\mathrm{Z}(s)\bigr)$ is a trivial invertible
sheaf, and there is an isomorphism $\mathcal{L}\simeq
\mathcal{O}\bigl(\mathrm{Z}(s)\bigr)$.

\begin{remark}
Generally, on  a regular stack not all invertible sheaves come
from divisors; those that do are precisely those possessing a
rational section that does vanish on open dense substack. For
example, if $G$ is a finite group, the group of divisors on the
associated classifying stack
$\mathcal{B}_{\mathbb{C}}G$ is trivial, while the Picard group is
the group of characters $G \to \mathbb{C}^*$ of $G$. In this case
a rational section is an invariant, and only the trivial
character has nonzero invariants.
\end{remark}

Now let us proceed with the proof of the lemma; set
$\mathcal{X}_K = \spec K \times_{\spec R} \mathcal{X}$.

Let us show that the restriction homomorphism $\pic \mathcal{X}
\to \pic \mathcal{X}_K$ is injective. Let $\mathcal{L}$ be an
invertible sheaf on $\mathcal{X}$ whose restriction to
$\mathcal{X}_K$ is trivial. Choose a nowhere vanishing section of
the restriction of $\mathcal{L}$ to $\mathcal{X}_K$; this will
extend to a nowhere vanishing section $s$ of $\mathcal{L}$ over
some open substack $\mathcal{U}$ of $\mathcal{X}$ containing the
fiber at infinity. Let $\mathcal{D}$ be the divisor on
$\mathcal{X}$ defined by $s$; then as we have seen $\mathcal{L}$
is isomorphic to $\mathcal{O}(\mathcal{D})$.
The support of $\mathcal{D}$ will be contained in a union of
closed fibers of the morphism $\mathcal{X} \to \spec R$; since
these fibers are integral we see that $\mathcal{D}$ is the
pullback of a divisor on $\spec R$, so $\mathcal{L}$ is the
pullback of an invertible sheaf on $\spec R$. But such an
invertible sheaf is alway trivial, because $R$ is a unique
factorization domain.

To prove surjectivity, take an invertible sheaf $\mathcal{M}$
over $\mathcal{X}_K$, and consider the quasicoherent sheaf $j_*
\mathcal{M}$ on $\mathcal{X}$, where $j \colon \mathcal{X}_K \to
\mathcal{X}$ is the natural morphism. We claim that the natural
homomorphism $j^*j_* \mathcal{M} \to \mathcal{M}$ is an
isomorphism. In fact, this is a local question in the smooth
topology of $\mathcal{X}$, so we may assume that $\mathcal{X}$ is
the spectrum of an $R$-algebra $A$; and then this follows from the
fact that $\mathcal{X}_K$ is the spectrum  of a localization of
$A$.

It follows from \cite[Proposition~15.4]{lm} that there exists a
coherent subsheaf $\mathcal{F}$ of $j_* \mathcal{M}$ whose
restriction to $\mathcal{X}_K$ coincides with $\mathcal{M}$. Then
the double dual
$\mathcal{F}^{\vee\vee}$ is a reflexive sheaf of rank~$1$ on a
regular stack, so it is invertible, and its restriction to
$\mathcal{X}_K$ is isomorphic to $\mathcal{M}$. This completes
the proof of the lemma.
\end{proof}

So, assume that $R$ equals a field $k$. From the
description of $\covsm$ in
Corollary~\ref{cor:describe-covsm} and from
\cite[Proposition~18]{Int} it  follows that
$\pic(\covsm)$ is equal to
$\chow^{1}_{\GL_{n+1}/\mmu[d]}\bigl(\formssm\bigr)$, the
codimension~1 component of the integral
$\GL_{n+1}/\mmu[d]$-equivariant Chow ring of $\formssm$.

Suppose that $G$ is an algebraic group over a field $k$, $V$ an
$l$-dimensional  representation of $G$, $X$ an open invariant subscheme of
$V$. If follows from \cite{Int} that the pullback
$\chow^1_G \eqdef \chow^1_G(\spec k) \to \chow^1_G(V)$ is an
isomorphism. Indeed, $\chow^1_G(\spec k) \simeq \chow^G_{-1}(\spec k)$ and
$\chow^1_G(V)\simeq \chow^G_{l-1}(V\times_{\spec k}\spec k)$ by
\cite[Proposition 4]{Int}; by \cite[Theorem 1]{Int} we get
$\chow^G_{-1}(\spec k)\simeq \pic^G(\spec k)$ and analogously for
$\chow^G_{l-1}(V\times_{\spec k}\spec k)$. Finally, by \cite[Lemma 2]{Int},
if $\pi : V\times_{\spec k}\spec k \to \spec k $ is the second projection,
then $\pi^{*}: \pic^G(\spec k) \to \pic^G( V\times_{\spec k}\spec k)$ is an
isomorphism and this yields the claim.

Again, $\chow_G^1(\spec k)$ is the equivariant Picard
group for the trivial action of $G$ over $\spec k$, that is, is
the group of characters $\widehat{G}$. Call $n$ the dimension
of $V$. From the usual exact sequence
    \[
    \chow^G_{n-1}(V \setminus X) \longrightarrow
    \chow^1_G(V) \longrightarrow \chow^1_G(X) \longrightarrow 0
    \]
we see that $\chow^1_G(X)$ is the quotient of $\widehat{G}$ by
the subgroup generated by the classes of the components of $V
\setminus X$ in codimension~1. In our case, the group of
characters $\widehat{\GL_{n+1}/\mmu[d]}$ is infinite
cyclic, while the locus $\Delta$ of singular forms is well
known to be irreducible, so $\chow^1_G(\formssm)$ is a cyclic
group, of order equal the the index of the subgroup generated by
the class of $\Delta$ in $\chow^1_{}\bigr(\forms\bigl) =
\widehat{\GL_{n+1}/\mmu[d]}$. To compute this index, first of all
notice that $\widehat{\GL_{n+1}/\mmu[d]}$ injects inside
$\widehat{\GL}_{n+1}$, that is generated by the determinant $\det
\colon \GL_{n+1} \to \gm$; since the intersection of $\mmu[d]$
with the kernel of the determinant has order $\gcd(d, n+1)$ it
follows that the index of $\widehat{\GL_{n+1}/\mmu[d]}$ inside
$\widehat{\GL}_{n+1}$ is $d/\gcd(d, n+1)$. In turn, if $\gm \into
\GL_{n+1}$ is the usual embedding, $\widehat{\GL}_{n+1}$ has
index $n+1$ in $\widehat{\mathbb{G}}_{\mathrm{m}}$; the
composite homomorphism $\gm \to \GL_{n+1}/\mmu[d]$ induces an
embedding $\widehat{\GL_{n+1}/\mmu[d]} \into
\widehat{\mathbb{G}}_{\mathrm{m}}$ of infinite cyclic groups with
index $(n+1)d/\gcd(d, n+1)$. The resulting action of $\gm$ on
$\forms$ is defined by the formula $\alpha\cdot f(x) =
f(\alpha^{-1}x) = \alpha^{-rd}f(x)$; thus the index of the
subgroup generated by the class of $\Delta$ in
$\widehat{\GL_{n+1}/\mmu[d]}$ equals the index of the subgroup of
the class of $\Delta$ in $\widehat{\mathbb{G}}_{\mathrm{m}}$ for
the action described above, multiplied by the rational number
$\gcd(d, n+1)/(n+1)d$.

Now, the action of $\gm$ described above is induced by the
standard action of $\gm$ defined by the usual formula
$\alpha\cdot f(x) = \alpha f(x)$ via the morphism $\gm \to \gm$
defined by $\alpha \mapsto \alpha^{-rd}$; hence the index of the
subgroup generated by the class of $\Delta$ in
$\widehat{\mathbb{G}}_{\mathrm{m}}$ for the action above is $rd$
times the class of $\Delta$ in
$\widehat{\mathbb{G}}_{\mathrm{m}}$ for the standard action. But
the class of $\Delta$ in $\widehat{\mathbb{G}}_{\mathrm{m}}$ for
the standard action is the degree of $\Delta$. Putting all this
together we get the following.

\begin{lemma}
If $k$ is a field, the Picard group of the stack $\covsm_k$
is cyclic of order equal to the degree of the hypersurface
$\Delta$ in $\forms_k$ consisting of singular forms, multiplied
by $r\gcd(d, n+1)/(n+1)$.
\end{lemma}

The hypersurface $\Delta$ in $\forms$ is well known to be defined
by a polynomial of degree $(n+1)(rd-1)^n$ (see for instance
\cite{GKZ}); the result would follow if we showed that this
polynomial is irreducible when the characteristic of $k$ does not
divide $2rd$.

Since $\Delta$ is a cone, we can compute its degree as the
degree of its projectivization $\overline{\Delta} \subseteq
\pforms$ (recall that $\pforms$ is the projective space of lines
in $\forms$). Call $N$ the dimension of $\pforms$. Let us
represent a point of $\mathbb{P}^n
\times \pforms$ as a pair $(x, F)$, and let us denote by $D$
the subscheme of $\mathbb{P}^n
\times \pforms$ defined by the homogeoneous
equations $\partial F/ \partial x_i = 0$ with $i = 0$,
\dots,~$i=n$; these are $n+1$ equations of bidegree $(rd-1, 1)$.
The projection $D \to \mathbb{P}^n$ makes $D$ into a
$\mathbb{P}^{N-n-1}$ bundle onto $\mathbb{P}^n$, hence $D$ is
smooth of codimension $n+1$, and a complete intersection. Call
$\xi$ and $\eta$ the classes in $\chow^1 (\mathbb{P}^n
\times \pforms)$ obtained by pulling back a hyperplane from
$\mathbb{P}^n$ and from $\pforms$ respectively; then the class
of $D$ in the Chow ring of $\mathbb{P}^n \times \pforms$ is
$\bigl((rd-1)\xi + \eta\bigr)^{n+1}$; a straightforward
calculation, applying projection formula reveals that its
pushforward to $\chow^1(\pforms)$
has degree $(n+1) (rd-1)^n$.

Because of Euler's formula, and because the degree $rd$ of a form
in $\forms$ is not divisible by the characteristic of $k$, if $(x,
F)$ is a point of $D$ then $x$ is a singular point of the
hypersurface defined by $F$; hence the image of $D$ in $\pforms$
is the projectivization $\overline{\Delta}$ of $\Delta$; hence to
conclude the proof is enough to show that $D$ is birational
onto $\overline{\Delta}$. Call $D_0$ the inverse
image of $D$ in $\formsnoto$; it is enough to show that
$D_0$ is birational onto its image in $\formsnoto$. We may
also assume that the base field $k$ is infinite. It is enough to
show that there exists a polynomial $F$ in $\formsnoto(k)$, whose
inverse image in $D_0$ is a single rational point
with the reduced scheme structure. Because of the definition
of $F$, this is equivalent to saying that $F$ has a single
singular point $p \in \mathbb{P}^n(k)$, and the ideal generated
by the partial derivates $\partial F/ \partial x_i$ is the
homogeneous ideal of $p$.

Take a polynomial $f \in k[x]$ in one variable of degree $rd$,
that has a double root in $0$ and no other multiple root. We set
$f = \sum_{i=1}^n a_i f(x_i)$, where $a_1$, \dots,~$a_n$ are
generic elements of $k$, and we call $F$ the homogeneous
polynomial of degree $rd$ whose dehomogenization is $f$. We
claim that $F$ has the desired property. One checks immediately
that $F$ has no singularity along the hyperplane at infinity.
Furthermore, if $\xi = (\xi_1, \dots, \xi_n)$ is a singular point
of $f$, then $0 = \partial f / \partial x_i(\xi) = a_i f'(\xi_i)$,
so $\xi_i$ is a zero of the derivative $f'$ of $f$. There are
only finitely many such zeros; hence, since the $a_i$ are
generic $f$ will not vanish on any $n$-tuple $(\xi_1, \dots,
\xi_n)$ where each $\xi_i$ is a zero of $f'$, and at least one of
them is different from $0$. So the only singularity of $f$ is at
the origin. But $\partial f / \partial x_i$ has the form $c_i
x_i + \text{higher order terms}$ with all the $c_i$ different
from $0$, so the partial derivatives
$\partial f / \partial x_i$ generate the ideal $(x_1, \dots,
x_n) \subseteq k[x_1, \dots , x_n]$. Again by Euler's formula
this implies that the ideal generated by the partial derivatives
$\partial F / \partial x_i$ for $i = 0$, \dots,~$n$ is the
homogeneous ideal $(x_1 - x_0, \dots, x_n - x_0)$, and this
completes the proof of Theorem~\ref{thm:picard}.
\end{proof}

\begin{remark} In particular, this says that the Picard group of
the stack of hyperelliptic curves $\covsm[1,2,g+1]_k$ over a
field $k$ of characteristic not dividing $2$ or $g+1$ is cyclic
of order $2(2g+1)$ if $g$ is even, and
$4(2g+1)$ if $g$ is odd.

When $g = 1$ one gets that that the Picard group of
$\covsm[1,2,1]_k$ is cyclic of order $12$; this immediately
reminds one of the famous result of Mumford in \cite{mum}, that the Picard
group of the stack $\mathcal{M}_{1,1}$ is cyclic of order $12$.
However, as we observe in Example \refall{ex:examples}{1},
$\covsm[1,2,1]_k$ is not isomorphic to $\mathcal{M}_{1,1}$. There
is a canonical morphism $\mathcal{M}_{1,1} \to \covsm[1,2,1]$,
sending a family
$\pi \colon E \to S$ to the uniform covering $E \to
\mathbb{P}\bigl(\pi_*
\mathcal{O}_E(2\Sigma) \bigr)$, where $\Sigma$ is the image of the
given section $S \to E$; a generator of $\mmu[2]$ acts like the
involution $e
\mapsto -e$ on $E$. This morphism induced a factorization $\mathcal{M}_{1,1}
\to  \covsm[1,2,1]_k \to \mathcal{M}_1$ of the morphism $\mathcal{M}_{1,1}
\to \mathcal{M}_1$ forgetting the section.

We claim that this morphism, although it is not an isomorphism,
induces an isomorphism of Picard groups. This can be seen as
follows. The Picard group of $\mathcal{M}_{1,1}$ is generated by
the first Chern class of the Hodge bundle on
$\mathcal{M}_{1,1}$. The Hodge bundle is already defined on the
stack
$\mathcal{M}_1$ of unpointed curves of genus~1, and the morphism
$\mathcal{M}_{1,1} \to \mathcal{M}_1$ forgetting the section
factors through $\covsm[1,2,1]_k$. Hence there is an element of
the Picard group of $\covsm[1,2,1]_k$ mapping into a generator of
the Picard group of
$\mathcal{M}_{1,1}$. Since both groups are cyclic of the same
order it follows that the pullback homomorphism is in fact an
isomorphism.
\end{remark}

\section{Cyclic triple coverings of $\mathbb{P}^1$}
\label{sec:triple}

In this section we study the stack of cyclic triple covers of the
projective line, with particular regard to the smooth ones.
General triple covers have been extensively studied starting from
\cite{M}.

\begin{definition}\call{def:triple-cyclic}
A \emph{cyclic triple cover} over a scheme $S$ consists of a
morphism of
$S$-schemes $f \colon X \to P$, together with an action of
$\mmu[3]$ over
$X$ leaving $f$ invariant, such that the following conditions are
satisfied.

\begin{enumeratea}

\item $P \to S$ is a conic bundle.

\itemref{triple1} The morphism $f$ is flat and finite, and induces an
isomorphism
$X/\mmu[3] \simeq P$.

\itemref{triple2} There exists an open subscheme $V \subseteq P$,
which intersects every fiber of $f \colon X \to P$, such that the
restriction
$f^{-1}(V) \to V$ is a $\mmu[3]$-torsor.

\end{enumeratea}
\end{definition}

Cyclic triple covers can be described by using an eigensheaf
decomposition, as for uniform cyclic covers. Consider the action
of
$\mmu[3]$ on the locally free sheaf $f_* \mathcal{O}_X$ on
$\mathcal{O}_P$; this will split as a sum of locally free sheaves
of
$\mathcal{O}_P$-modules $\mathcal{L}_0 \oplus \mathcal{L}_1 \oplus
\mathcal{L}_2$, where $\mathcal{L}_i$ is the subsheaf of $f_*
\mathcal{O}_X$ of sections $s$ such that the action of $\mmu[3]$
can be describes as $(t,s) \mapsto t^i s$.
Condition~\refpart{def:triple-cyclic}{triple1} of the definition
insures that $\mathcal{L}_0 = \mathcal{O}_P$, while flatness and
condition~\refpart{def:triple-cyclic}{triple2} imply that
$\mathcal{L}_1$ and $\mathcal{L}_2$ are invertible sheaves. The
algebra structure on $f_* \mathcal{O}_X$ induces homomorphisms of
sheaves of $\mathcal{O}_P$-algebras
     \[
     \phi_1\colon \mathcal{L}_1^{\otimes 2} \longrightarrow
     \mathcal{L}_2,
    \quad
     \phi_2\colon \mathcal{L}_2^{\otimes 2} \longrightarrow
     \mathcal{L}_1 \quad
     \text{and} \quad
     \phi_{12} \colon \mathcal{L}_1 \otimes \mathcal{L}_2
     \to \mathcal{O}_P
     \]
that are injective on every fiber of $P \to S$. These
homomorphism determine the algebra structure completely; the
covering is uniform if and only if
$\phi_1$ is an isomorphism. The algebra structure also gives
homomorphisms
$\mathcal{L}_1^{\otimes 3}
\to
\mathcal{O}_P$ and
$\mathcal{L}_2 \to \mathcal{O}_P$ which are injective on every
fiber of $P
\to S$: this says that the degrees of the
$\mathcal{L}_i$ on each fiber of $P \to S$ can not be positive.
We will assume that these degrees are constant on $S$, and we
will call their opposites $d_1$ and $d_2$ the
\emph{branch degrees} of the triple covering. These
branch degrees are subject to the obvious constraints $0
\le d_1 \le 2d_2$, $0
\le d_2 \le 2d_1$.

The stack of cyclic triple covers with branch degrees $d_1$
and
$d_2$ will be denoted by $\triple$; we have $\cov[1,3,d] =
\triple[1,3;d, 2d]$. We will denote by $\triplesm$ the full
subcategory of $\triple$ whose objects are triple cyclic covers
$X \to P \to S$ such that $X$ is smooth over $S$.

Of course all the definition above generalize to higher
dimension, and we could consider categories of cyclic triple
covers of $\mathbb{P}^n$ for any $n$; the main reason why we do
not do this is that such a cover will never be smooth, unless $n
= 1$, or the cover becomes uniform after twisting the action by an
automorphism of $\mmu[3]$ (see Remark~\ref{rmk:why-not}).

Here is an alternate description of $\triple$. We call
$\alttriple$ the category whose objects are quintuples $(P \to S,
\mathcal{L}_1,
\mathcal{L}_2, \phi_1, \phi_2)$, where $P \to S$ is a
Brauer--Severi scheme of rank~1, $\mathcal{L}_1$ and
$\mathcal{L}_2$ are invertible sheaves on $P$, whose degrees on
each fiber of $P \to S$ are $-d_1$ and
$-d_2$ respectively, while $\phi_1\colon \mathcal{L}_1^{\otimes 2}
\to
\mathcal{L}_2$ and $\phi_2\colon \mathcal{L}_2^{\otimes 2} \to
\mathcal{L}_1$ are homomorphism of sheaves of
$\mathcal{O}_P$-modules that are injective on all the fibers of
$P \to S$. The arrows are defined in the obvious way.

The construction above yields a functor $\triple \to \alttriple$;
We claim that this is an equivalence of fibered categories over
the category of schemes. This is an easy consequence of the
following.

\begin{lemma}\label{lem:unique-extension}
Let $Y$ be a scheme, $\mathcal{L}_1$,
$\mathcal{L}_2$ invertible sheaves on $Y$, with homomorphisms
$\phi_1\colon \mathcal{L}_1^{\otimes 2} \to
\mathcal{L}_2$ and $\phi_2\colon \mathcal{L}_2^{\otimes 2} \to
\mathcal{L}_1$. Then
$\phi_1$ and $\phi_2$ extend to  a unique structure of
associative and commutative $\mathcal{O}_Y$-algebra on the
$\mathcal{O}_Y$-sheaf
$\mathcal{O}_Y \oplus \mathcal{L}_1 \oplus \mathcal{L}_2$.
\end{lemma}

\begin{proof}
This is a local statement in the Zariski topology, so we may
assume that $\mathcal{L}_1$ and $\mathcal{L}_2$ have global
generators $t_1$ and $t_2$. The homomorphisms $\phi_1$ and
$\phi_2$ correspond to two section $f_1$ and $f_2$ of
$\mathcal{O}_Y$ with $\phi_1(t_1 \otimes t_1) = f_1t_2$ and
$\phi_2(t_2 \otimes t_2) = f_2t_1$. Set $\mathcal{A} \eqdef
\mathcal{O}_Y \oplus \mathcal{L}_1 \oplus \mathcal{L}_2$; to
extend $\phi_1$ and
$\phi_2$ to a bilinear symmetric product $\mathcal{A}
\otimes_{\mathcal{O}_Y} \mathcal{A} \to \mathcal{A}$ with
identity $1$ we need to add the data of a homomorphism
$\mathcal{L}_1 \otimes \mathcal{L}_2 \to \mathcal{O}_Y$,
corresponding to a third section $h$ of $\mathcal{O}$ (the image
of $t_1 \otimes t_2$). Then a lengthy but straightforward
calculation reveals that this product is associative if and only
if $h = f_1f_2$, and this clearly implies the result.
\end{proof}

A cyclic triple cover $X \to P$ has two associated branch
divisors in $P$, given by the two homomorphisms
$\phi_1 \colon \mathcal{L}_1^{\otimes 2} \to \mathcal{L}_2$ and
$\phi_2\colon \mathcal{L}_2^{\otimes 2} \to \mathcal{L}_1$, whose
degrees are respectively $2d_1 - d_2$ and $2d_2 - d_1$. We say
that the triple cover is \emph{smooth} if $X$ is smooth over $S$;
we denote by $\triplesm$ the open substack of $\triple$ whose
objects are smooth triple covers.

Since $X$ is smooth over $S$, to check smoothness it is enough to
check that the geometric fibers are smooth.

\begin{proposition}\label{prop:trip-smooth}
A cyclic triple cover $X \to P$ over a field is smooth if and only
its two branch divisors have no multiple points, and are disjoint.
\end{proposition}

\begin{proof}
This follows from Proposition 3.1 of \cite{Cov}.

One can also proceed as follows. Choose an
open subset
$U$ of $P$ with nonvanishing section $t_1$ and $t_2$ of
$\mathcal{L}_1$ and $\mathcal{L}_2$ respectively. The
$\mathcal{O}_P$ algebra $\mathcal{O}_X$ is defined over $U$ by
the equations $t_1^2 = f_1t_2$, $t_2^2 = f_2t_1$,
$t_1t_2 = f_1f_2$ (see the proof of
Lemma~\ref{lem:unique-extension}). A straightforward calculation
using the Jacobian criterion proves that $X$ is smooth over $S$
if and only if $f_1$ and $f_2$ have no multiple zero and no
common zero.
\end{proof}

\begin{remark}\label{rmk:why-not}
One could build a similar theory for projective spaces of
dimension higher than 1; then a similar argument would show that
a triple cover $X \to P$ is smooth if and only if its two
branch divisors are smooth and do not intersect. However,
in rank greater than 1, this would mean that one of the two
divisors must be empty, so that either $d_2 = 2 d_1$, and the
triple cover is in fact uniform, or $d_1 = 2d_2$, and the triple
cover becomes uniform after twisting the action of $\mmu[3]$ by
the nontrivial automorphism of $\mmu[3]$. See
Remark~\ref{rmk:why-not1}.
\end{remark}

Using this description of $\triple$ we can prove the following.
Consider the embedding $\mmu[d_1] \times \mmu[d_2] \subseteq \gm
\times \GL_2$ as a normal subgroup scheme given by
     \[
     (\alpha_1, \alpha_2) \mapsto (\alpha_2/ \alpha_1, \alpha_1
     \mathrm{I}_2);
     \]
call $\Gamma(d_1, d_2)$ the quotient.

\begin{theorem}\label{thm:quotient-triple}
$\triple$ is isomorphic to the
quotient stack
     \[
     [ \formsnoto[1, 2d_1 - d_2] \times \formsnoto[1, 2d_2 - d_1]/
     \Gamma(d_1, d_2)]
     \]
by the action given by the formula
     \[
     [\alpha, A]\cdot \bigl(f_1 (x), f_2 (x)\bigr) =
     \bigl(\alpha^{d_2}f_1(A^{-1}x),
    \alpha^{-2d_2}f_2(A^{-1}x)\bigr).
     \]
Furthermore, $U$ is the open subscheme of $\formsnoto[1,3,2d_1 -
d_2]
\times \formsnoto[1,3,2d_2 - d_1]$ consisting of pairs of forms
without multiple roots and no common root, then $\triplesm$ is
isomorphic to the quotient $[U/\Gamma(d_1, d_2)]$.
\end{theorem}

\begin{proof}
We follow closely the strategy of the
proof of Theorem~\ref{thm:quotient}, and use the alternative
description of $\triple$ given by $\alttriple$.

Consider the auxiliary fibered category $\triplee$: if $S$ is a
scheme, an object of
$\triplee(S)$ is a quintuple $(P \to
S,
\mathcal{L}_1,
\mathcal{L}_2, \phi_1, \phi_2)$ giving an object of
$\alttriple(S)$, plus the choice of two isomorphisms
over $S$
    \[
    \lambda_1 \colon
    \bigl(\mathbb{P}^1_S,\mathcal{O}(-d_1)\bigr)
    \simeq (P,\mathcal{L}_1)
    \quad\text{and} \quad
    \lambda_2 \colon
    \bigl(\mathbb{P}^1_S,\mathcal{O}(-d_2)\bigr)
    \simeq  (P,\mathcal{L}_2)
    \]
such that the
restriction of
$\lambda_2^{-1} \circ \lambda_1$ to $\mathbb{P}^1_S$ induces
the identity on
$\mathbb{P}^1_S$. The arrows are arrows in $\triple$ preserving
the two isomorphisms $\lambda_1$ and $\lambda_2$. Objects in
the fiber $\triplee(S)$ have no notrivial automorphism, so
$\triplee$ is equivalent to a functor.

Let $(P \to S, \mathcal{L}_1, \mathcal{L}_2, \phi_1, \phi_2,
\lambda_1, \lambda_2)$ be an object of $\triplee(S)$; consider
the composition
    \[
    \bigl(\mathbb{P}^1_S, \mathcal{O}(-2d_1)\bigr)
    \xrightarrow{\lambda_1^{\otimes 2}}
    \bigl(P, \mathcal{L}_1^{\otimes 2}\bigr)
    \xrightarrow{\phi_1}
    \bigl(P, \mathcal{L}_2\bigr)
    \xrightarrow{\lambda_2^{-1}}
    \bigl(\mathbb{P}_S^1, \mathcal{O}(-d_2)\bigr)
    \]
(if
$\lambda_1$ is given by a pair $(\mu_1, \rho_1)$, where
$\mu_1 \colon P \simeq \mathbb{P}^1_S$ is an isomorphism of
$S$-schemes and $\rho_1\colon \mathcal{O}(-1) \simeq
\mu_1^*\mathcal{O}(-d_1)$ is an isomorphism of sheaves of
$\mathcal{O}_{\mathbb{P}^1_S}$-modules, we write
$\lambda_1^{\otimes 2}$ for the pair $(\mu_1, \rho_1^{\otimes
2})$). This gives a homomorphism of sheaves
$\mathcal{O}(-2d_1) \to
\mathcal{O}(-d_2)$ lying over the identity of $\mathbb{P}^1_S$,
that does not vanish identically on any fiber, and therefore an
element of $\formsnoto[1, 2d_1 - d_2](S)$. Analogously we
construct an element of $\formsnoto[1, 2d_2 - d_1](S)$; this
defines a base-preserving functor
    \[
    \triplee \longrightarrow \formsnoto[1, 2d_1 - d_2]
    \times \formsnoto[1, 2d_2 - d_1].
    \]

To define an inverse, we send an object $(f_1, f_2)$ of
$\formsnoto[1, 2d_1 - d_2](S) \times \formsnoto[1, 2d_2 -
d_1](S)$ into the object
    \[
    \bigl(\mathbb{P}^1_S \to S, \mathcal{O}(-d_1),
    \mathcal{O}(-d_2), \phi_1, \phi_2, \id, \id \bigr)
    \]
where $\phi_1\colon \mathcal{O}(-2d_1) \to \mathcal{O}(-d_2)$ and
$\phi_2\colon \mathcal{O}(-2d_2) \to \mathcal{O}(-d_1)$ are given
by multiplication by $f_1$ and $f_2$ respectively. It is
straightforward to check that this gives a quasi-inverse to the
functor above; so $\triplee$ is equivalent to $\formsnoto[1, 2d_1
- d_2] \times \formsnoto[1, 2d_2 - d_1]$.

Consider the
functor $\biautomo{-d_1}{-d_2}$ from schemes to groups, sending
a scheme $S$ into the the group of automorphisms of the triple
$\bigl(\mathbb{P}^1_S, \mathcal{O}(-d_1), \mathcal{O}(-d_2)\bigr)$
over the identity on $S$. This is a sheaf in the fppf topology.
It can also be thought of as the fiber product
    \[
    \automo{-d_1} \times_{\underaut{\mathbb{P}^1_\mathbb{Z}}} \automo{-d_2}.
    \]
Since $\underaut{\mathbb{P}^1_\mathbb{Z}}$ is $\PGL_{2, \mathbb{Z}}$,
and, according to the discussion in the proof of
Theorem~\ref{thm:quotient}, $\automo{-d}$ is isomorphic to the
quotient $\GL_{2, \mathbb{Z}}/\mmu[d, \mathbb{Z}]$, we see that
we have an isomorphism of functors
    \[
    \biautomo{-d_1}{-d_2} \simeq
    \GL_{2} \times_{\PGL_2} \GL_2/\mmu[d_1] \times \mmu[d_2]
    \]
where $\mmu[d_1]$ is embedded diagonally in the first copy of
$\GL_2$, $\mmu[d_2]$ in the second. We also have an isomorphism
$\gm \times \GL_2 \simeq \GL_{2} \times_{\PGL_2} \GL_2$, where a
section $(\alpha, A)$ of $\gm \times \GL_2$ over some scheme is
sent into the pair $(A, \alpha A)$; the embedding $\mmu[d_1]
\times \mmu[d_2] \subseteq \GL_{2} \times_{\PGL_2} \GL_2$ gives
an embedding $\mmu[d_1] \times \mmu[d_2] \subseteq \gm \times
\GL_2$ given by the formula
    \[
    (\alpha_1, \alpha_2) \mapsto
    (\alpha_2/ \alpha_1, \alpha_1\mathrm{I});
    \]
in this way we obtain an isomorphism of $\biautomo{-d_1}{-d_2}$
with $\Gamma(d_1, d_2)$.

There is a left action of $\biautomo{-d_1}{-d_2}$ on $\triplee$;
if $(P \to S, \mathcal{L}_1, \mathcal{L}_2, \phi_1, \phi_2,
\lambda_1, \lambda_2)$ is an object of $\triplee (S)$ and
$(\alpha_1, \alpha_2)$ is an object of $\biautomo{-d_1}{-d_2}
(S)$, we set
    \begin{align*}
    (\alpha_1, \alpha_2) \cdot
    &(P \to S, \mathcal{L}_1, \mathcal{L}_2, \phi_1, \phi_2,
    \lambda_1, \lambda_2)\\
    = {}&(P \to S, \mathcal{L}_1, \mathcal{L}_2,
    \phi_1,
    \phi_2,
    \lambda_1 \circ \alpha_1^{-1},
    \lambda_2 \circ \alpha_2^{-1}).
    \end{align*}
Furthermore, given two invertible sheaves
$\mathcal{L}_1$ and $\mathcal{L}_2$ on
$P \to S$ with degrees $-d_1$ and $-d_2$ on every geometric fiber,
there is an fppf covering $S' \to S$, such that the pullback of the triple
$(P,\mathcal{L}_1,\mathcal{L}_2)$ to $S'$ is isomorphic to
$\bigl(\mathbb{P}^1_{S'},\mathcal{O}(-d_1), \mathcal{O}(-d_2)\bigr)$;
this fact, plus
descent theory, implies that the forgetful morphism
$\triplee \to \triple$ makes $\triplee$ into a principal bundle
with structure group
$\biautomo{-d_1}{-d_2}=\Gamma(d_1, d_2)$.

The action of $\Gamma(d_1, d_2)$ on $\triplee$ gives an action of the
structure group $\biautomo{-d_1}{-d_2}$ on $\formsnoto[1, 2d_1 - d_2] \times
\formsnoto[1, 2d_2 - d_1]$, via the equivalence above; hence
$\triple$ is equivalent to the quotient stack
    \[
    [\formsnoto[1, 2d_1 - d_2] \times \formsnoto[1, 2d_2 - d_1]/
    \Gamma(d_1, d_2)];
    \]
we have only left to write this action explicitly.

First of all, restrict the attention to the
first component
$\formsnoto[1, 2d_1 - d_2]$ and consider the action of $\GL_2
\times_{\PGL_2}
\GL_2$ on $\H^0(\mathcal{O}\bigl(2d_1-d_2)\bigr)$. Fix a section
$h_1$ in
$\H^0\bigl(\mathbb{P}^1_{\mathbb{Z}}, \mathcal{O}(d_2)\bigr)$ that
does not vanish on any fiber of $\mathbb{P}^1_{\mathbb{Z}} \to
\spec\mathbb{Z}$: any element
$f_1\in
\H^0\bigl(\mathcal{O}(2d_1 -d_2)\bigr)$ can be written uniquely as
$f_1(x)=g_1(x)/h_1(x)$, where $g_1\in \H^0(
\mathcal{O}(2d_1))$.

Since the multiplication map $\H^0(\mathcal{O}(2d_1-d_2))\times
\H^0(\mathcal{O}(d_2)) \to \H^0(\mathcal{O}(2d_1))$ is $\GL_2
\times_{\PGL_2}\GL_2$-equivariant, it follows immediately that
the action of the group
$\GL_2 \times_{\PGL_2}\GL_2$ on $H^0(\mathcal{O}(2d_1-d_2))$ is
described by the formula
    \[
    (A_1,A_2)f_1(x) = g_1(A_1^{-1}x)/h_1(A_2^{-1}x).
    \]
Under the isomorphism $\gm \times \GL_2 \to
\GL_2 \times_{\PGL_2} \GL_2 $, given by $(\alpha,A) \mapsto (A,
\alpha A)$, the action of $\gm \times \GL_2$ on
$H^0(\mathcal{O}(2d_1-d_2))$ can be  written as
    \[
    (\alpha, A)f_1(x) = g_1(A^{-1}x)/h_1(\alpha^{-1}A^{-1}x)
    =\alpha^{d_2}f_1(A^{-1}(x)).
    \]
It is easy to check that this action descends to the
quotient group $\Gamma(d_1,d_2)$. Analogously, one
shows that the action  of
$\Gamma(d_1,d_2)$ on the second component
$\formsnoto[1, 2d_2 - d_1]$ is  given by the formula
\[
(\alpha,
A)f_2(x) = \alpha^{-2d_2}f_2(A^{-1}(x)).
\]
This completes the proof of the first statement.

The last statement is an easy consequence of the
definition of the stack of smooth triple covers $\triplesm$ and of Proposition
\ref{prop:trip-smooth}.
\end{proof}

In particular, the stack $\triple$ is a smooth irreducible Artin
stack over
$\spec\mathbb{Z}$ of dimension $(2d_1 - d_2 + 1)(2d_2 - d_1 + 1)
- 5$.

Now we give a presentation of the Picard group of the stack
$\triplesm$.

\begin{theorem}
Assume that $d_1$ and $d_2$ are positive. Let R be a unique
factorization domain, such that the characteristic of its
quotient field does not divide
$2(2d_1-d_2)(2d_2-d_1)$. Then the Picard group
$\pic\bigr(\triplesm_R\bigl)$ is a group with two generators
$v_1$  and $v_2$ and three relations.
\begin{enumeratea}

\item If $d_1$ is odd, the three relations are:
\begin{gather*}
(2d_1-d_2-1)(2v_1-(d_2+2)v_2),\\
(2d_2-d_1-1)(4v_1 - (2d_2+1)v_2),\\
(-5d_1 + 4d_2)v_1 + \frac{4d_1 - 5d_2(d_1+1) -
4d_2^2}{2}v_2.
\end{gather*}

\item If $d_1$ and $d_2$ are both even, the three relations are
\begin{gather*}
2(2d_1-d_2-1)(v_1 - 2v_2),\\
2(2d_2-d_1-1)(2v_1 - v_2),\\
(4d_2-5d_1)v_1+(4d_2 - 5d_2)v_2.
\end{gather*}

\end{enumeratea}
\end{theorem}

\begin{remark}
As one sees immediately, twisting the action by the nontrivial
automorphism $\mmu[3] \simeq \mmu[3]$ gives a canonical
isomorphism of stacks $\triple \simeq
\triple[1,3;d_2,d_1]$; hence the theorem above describes the
Picard group of $\triple$ even when $d_1$ is even and $d_2$ is
odd.
\end{remark}

\begin{proof}
The proof is very similar to the proof of
Theorem~\ref{thm:picard}, with some added complications. Again,
Lemma~\ref{lem:Picard-to-generic} allows us to reduce to the case
that $R$ is a field.

The Picard group $\pic(\triplesm)$ is isomorphic to the
codimension~1 component
$\chow^{1}_{\Gamma(d_1,d_2)}(U)$ of the equivariant Chow ring of
the open subscheme $U$ of $\formsnoto[1, 2d_1
- d_2]\times \formsnoto[1, 2d_2 - d_1]$ consisting of the
complement of three hypersurfaces $\Delta_1$, $\Delta_2$ and $Z$:
the first two are the inverse images of the discriminant
hypersurfaces of $\forms[1, 2d_1 - d_2]$ and $\forms[1, 2d_2 -
d_1]$ respectively, while the geometric points of the third
consists of pairs of forms with a common zero (again, this 
description of $U$ comes from Proposition~\ref{prop:trip-smooth}). As 
in the proof of
Theorem~\ref{thm:picard}, this means that
$\chow^{1}_{\Gamma(d_1,d_2)}(U)$ is the quotient of
    \begin{align*}
    \chow^{1}_{\Gamma(d_1,d_2)}\bigl(\formsnoto[1, 2d_1 -
    d_2]\times \formsnoto[1, 2d_2 - d_1]\bigr) &=
    \chow^{1}_{\Gamma(d_1,d_2)}(\spec k)\\
    &= \widehat{\Gamma(d_1,d_2)}
    \end{align*}
by the subgroups generated by the classes of the three
hypersurfaces. The character group $ \widehat{\Gamma(d_1,d_2)}$
is the kernel of the restriction homomorphism
    \[
    \widehat{\gm \times \GL_2}
    \longrightarrow \widehat{\mmu[d_1] \times \mmu[d_2]};
    \]
The character group of $\gm \times \GL_2$ is generated by the
projection $e_1 \colon \gm \times \GL_2 \to \gm$ and by the
homomorphism $e_2\colon  \gm \times \GL_2 \to \gm$ defined by
$(\alpha, A) \mapsto \det A$. If we denote by $\epsilon_1$ and
$\epsilon_2$ the generators of $\widehat{\Gamma(d_1,d_2)}$
corresponding to the projection onto $\mmu[d_1]$ and $\mmu[d_2]$
followed by the embedding into $\gm$, the restriction
homomorphism sends $e_1$ into $\epsilon_2 - \epsilon_1$, and
$e_2$ into $2\epsilon_1$; from this we see that the kernel of the
restriction homomorphism is the subgroup of elements of
$\widehat{\gm \times \GL_2}$ of the form $x_1e_1 + x_2e_2$, where
$x_1$ and $x_2$ are integers with $x_1 \equiv 2x_2 \pmod{d_1}$
and $x_1 \equiv 0 \pmod{d_2}$. If $d_1$ is odd, then a basis for
the kernel is given by
    \[
    v_1 = d_2e_1 + \frac{(d_1 + 1)d_2}{2}e_2 \qquad
    \text{and} \quad
    v_2 = d_1e_2,
    \]
while if $d_1$ and $d_2$ are both even a basis is
    \[
    v_1 = d_2e_1 + \frac{d_2}{2}e_2 \qquad
    \text{and} \quad
    v_2 = \frac{d_1}{2}e_2.
    \]
So the Picard group of $\triple$ is generated by two elements
$v_1$ and $v_2$, with three relations, obtained by expressing the
classes of the three hypersurfaces as linear combinations of
$v_1$ and $v_2$.

To do this we use the following lemma. A \emph{cone} in
$\forms[2d_1 - d_2]_k \times \forms[2d_1 - d_2]_k$ is a closed
subscheme that is invariant under the actions of $\gm \times
\gm$ defined by $(t_1, t_2)\cdot(f_1, f_2) = (t_1
f_1, t_2 f_2)$. The integral cones
correspond to the integral subschemes of $\pforms[2d_1 - d_2]_k
\times
\pforms[2d_1 - d_2]_k$, and, as such, they have a bidegree.

\begin{lemma}\label{lem:formula-class}
Let $S$ be an integral cone of codimension~$1$ in $\forms[2d_1 -
d_2]_k \times \forms[2d_1 - d_2]_k$ of bidegree $(a_1, a_2)$ that
is  invariant under the action of $\Gamma(d_1, d_2)$.
\begin{enumeratea}

\item If $d_1$ is odd, the integer $4a_2 d_2 - a_2d_1 - a_1d_1d_2$
is divisible by $2d_1$, and the class of $S$ in
$\widehat{\Gamma(d_1, d_2)}$ is
    \[
    (a_1 - 2a_2)v_1 + (-a_1 + a_2d_2 + a_2a_1d_2/2)v_2
    \]

\item If $d_1$ and $d_2$ are both even, the integer
$4a_2d_2$ is divisible by $d_1$, and the class of $S$ in
$\widehat{\Gamma(d_1, d_2)}$ is
    \[
    (a_1 - 2 a_2)v_1 + (-2 a_1 + a_2)v_2.
    \]

\end{enumeratea}
\end{lemma}

\begin{proof}
Let $\Phi$ be a generator of the ideal of $S$. Saying that the
$S$ has bidegree $(a_1, a_2)$ is the same as saying that
$\Phi(t_1f_1, t_2f_2) = t_1^{a_1}t_2^{a_2}\Phi(f_1, f_2)$ for any
$(t_1, t_2)$ in $\gm \times \gm$ and any $(f_1, f_2)$ in
$\forms[2d_1 - d_2]_k \times \forms[2d_1 - d_2]_k$. On the other
hand, since $S$ is also invariant for the action of $\gm \times
\GL_2$, we must have a formula of the type
    \[
    \Phi\bigl(\alpha^{d_2}f_1(A^{-1}x),
    \alpha^{-2d_2}f_2(A^{-1}x)\bigr) =
    \alpha^{n_1} (\det A) ^{n_2} \Phi(f_1,f_2);
    \]
furthermore, in this case, the class of $S$ is $n_1e_1 + n_2e_2$, where
$e_1$ and $e_2$ are generators of $\widehat{\gm \times \GL_2}$.
To compute the integers $n_1$ and $n_2$ we set $A = \beta
\mathrm{I}_2$, where $\beta$ is scalar, so that $\det A =
\beta^2$. We get
    \begin{align*}
    \Phi\bigl(\alpha^{d_2}f_1(\beta^{-1}x),
    \alpha^{-2d_2}f_2(\beta^{-1}x)\bigr) &=
    \Phi\bigl(\alpha^{d_2}\beta^{-2d_1 + d_2}f_1(x),
    \alpha^{-2d_2}\beta^{-2d_1 + d_2}f_2(x)\bigr)\\
    &= \alpha^{(a_1 - 2a_2)d_2}
    \beta^{-a_1(2d_1 - d_2) - a_2(2d_2 - d_1)} \Phi(f_1, f_2)
    \end{align*}
hence
    \begin{multline*}
    \Phi\bigl(\alpha^{d_2}f_1(A^{-1}x),
    \alpha^{-2d_2}f_2(A^{-1}x)\bigr)\\ =
    \alpha^{(a_1 - 2a_2)d_2} (\det A) ^{-(a_1(2d_1 - d_2) +
    a_2(2d_2 - d_1))/2} \Phi(f_1,f_2)
    \end{multline*}
and from this we obtain that the class of $S$ in $\widehat{\gm \times
\GL_2} \simeq \mathbb{Z}^2$ is
    \[
    (a_1 - 2a_2)d_2 e_1 -
    \frac{a_1(2d_1 - d_2) + a_2(2d_2 - d_1)}{2}e_2.
    \]

The result follows by expressing this class as a linear
combination of $v_1$ and $v_2$.
\end{proof}

This reduces the problem to computing the
bidegrees of the three hypersurfaces. The hypersurface $\Delta_1$
is the pullback of the discriminant hypersurface from the first
factor
$\forms[2d_1-d_2]$, and we have seen in the proof of
Theorem~\ref{thm:picard} that this is integral with degree
$2(2d_1-d_2-1)$. Hence the bidegree of $\Delta_1$ is
$\bigl(2(2d_1-d_2-1), 0\bigr)$, and, if we plug this in the
formulas of Lemma~\ref{lem:formula-class} we obtain that the
class of
$\Delta_1$ is
    \[
    2(2d_1-d_2-1)v_1 - (d_2+2)(2d_1+d_2-1)v_2
    \]
when $d_1$ is odd, and
    \[
    2(2d_1-d_2-1) v_1 - 4(2d_1-d_2-1)v_2
    \]
when $d_1$ and $d_2$ are both even. This gives us our first
relation.

The second is obtained similarly, by setting $a_1 = 0$ and $a_2 =
2(2d_2 - d_1 -1)$ in the formulas; the result is
    \[
    2(2d_1 - d_2 - 1)v_1 - 4(2d_1 - d_2 - 1)v_2
    \]
when $d_1$ is odd, and
    \[
    -4(2d_2 - d_1 - 1)v_1 + (2d_2+1)(2d_2 - d_1 - 1)v_2
    \]
when $d_1$ and $d_2$ are both even.

To calculate the bidegree of $Z$, consider the subscheme
$\widetilde{Z}$ of $\mathbb{P}^1 \times
\pforms[2d_1-d_2] \times \pforms[2d_2-d_1]$ consisting of triples
$(p, f_1, f_2)$, where $p$ is a point in $\mathbb{P}^1$ and
$f_1$, $f_2$ are forms vanishing at $p$. Then $\widetilde{Z}$ is
a smooth subscheme of codimension~$2$, and the projection
    \[
    \mathbb{P}^1 \times \pforms[2d_1-d_2] \times \pforms[2d_2-d_1]
    \longrightarrow \pforms[2d_1-d_2] \times \pforms[2d_2-d_1]
    \]
maps $\widetilde{Z}$ birationally onto $Z$. If we denote by
$\eta$, $\xi_1$ and $\xi_2$ the pullbacks to $\mathbb{P}^1 \times
\pforms[2d_1-d_2] \times \pforms[2d_2-d_1]$ of the first Chern
classes of $\mathcal{O}(1)$ on $\mathbb{P}^1$,
$\pforms[2d_1-d_2]$ and $\pforms[2d_2-d_1]$ respectively, then
the class of $\widetilde{Z}$ in the Chow ring of  $\mathbb{P}^1
\times \pforms[2d_1-d_2] \times \pforms[2d_2-d_1]$ is
    \[
    \bigl((2d_1 - d_2)\eta + \xi_1\bigr)
    \bigl((2d_2 - d_1)\eta + \xi_2\bigr)
    \]
by pushing this class forward to $\pforms[2d_1-d_2] \times
\pforms[2d_2-d_1]$, using projection formula, we see that the 
bidegree of $Z$ is $(2d_2 -
d_1, 2d_1 - d_2)$. Again we use the formulas of
Lemma~\ref{lem:formula-class} to obtain the third relation
    \[
    (-5d_1 + 4d_2)v_1 + \frac{4d_1 - 5d_2(d_1+1) -
    4d_2^2}{2}v_2
    \]
when  $d_1$ is odd, and
    \[
    (-5 d_1 + 4d_2)v_2 + (4 d_1 - 5 d_2)v_2
    \]
when $d_1$ and  $d_2$ are both even.
\end{proof}

\begin{remark}\label{rmk:general-triple}

There are two possible generalization of this theory. Given
that general flat covers of $\mathbb{P}^1$ seem completely out
of reach, one could study the stack of cyclic covers of a conic of
degree $r$ for fixed $r$, or the stack of general triple covers.

For general triple covers of conics, one can use the description
of \cite{M}. This is work in progress of Marco Barone, a
student of one of us (Vistoli). There is one difficulty: although
every locally free sheaf of rank~2 on $\mathbb{P}^1_k$ is
isomorphic to
$\mathcal{O}(m) \oplus \mathcal{O}(n)$, if $k$ is a field, this is
not true over an arbitrary base. This means that we one can
mimick the construction of Theorem~\ref{thm:quotient-triple}, and
use the results of \cite{M} to study the stack of triple
coverings $f \colon X \to P$ where $P \to S$ is a conic bundle and
the kernel of the trace map $f_* \mathcal{O}_X \to \mathcal{O}_P$
is assumed to be locally isomorphic to $\mathcal{O}(m) \oplus
\mathcal{O}(n)$ for fixed $m$ and $n$; but removing this
unpleasant restriction requires a new idea.

For
general cyclic covers (or, more generally, covers that are
generically torsors under a finite diagonalizable group) there is
the theory created by Pardini (see \cite{Cov}). Using her
``reduced building data'' one can describe the stack of all
cyclic smooth covers of Brauer--Severi varieties as a quotient of
an open subsets of a representation of a quotient of a product of
general linear groups; but for non-smooth covers her description
does not work in general, and only yields a stack that is
birational to the stack of cyclic covers.
\end{remark}


\begin{thebibliography}{[G-K-MP98]}

\bibitem[Cat84]{C} F. Catanese: \emph{On the moduli spaces of
surfaces of general type}, J. Diff. Geom. {\bf 19} (1984),
483--515.

\bibitem[Ed-Gr98]{Int} D. Edidin and W. Graham:
\emph{Equivariant Intersection Theory}, Invent. Math. {\bf 131}
(1998), 595--634.

\bibitem[G-Z-K94]{GKZ} I.M. Gelfand, M.M. Kapranov, A.V.
Zelevinsky:
\emph{Discriminants, resultants and multidimensional
determinants}, Birkh\"auser Boston, (1994).

\bibitem[Kir85]{Kir} F. Kirwan: \emph{Partial desingularisations
of quotients of nonsingular varieties and their Betti numbers},
Ann. of Math. {\bf 122} (1985), 41--85.

\bibitem[L-MB00]{lm} G. Laumon, L. Moret-Bailly: \emph{Champs
alg\'ebrique}, Springer Verlag (2000).

\bibitem[Mir85]{M} R. Miranda: \emph{Triple covers in algebraic
geometry}, Amer. J. of Math. {\bf 107} (1985),113--1158.

\bibitem[Mum65]{mum} D. Mumford: \emph{Picard groups of moduli problems},
Arithmetical Algebraic Geometry (Proc. Conf. Purdue Univ., 1963), Harper \&
Row, New York (1965), 33--81.


\bibitem[Par91]{Cov} R. Pardini: \emph{Abelian covers of algebraic
varieties}, J. reine agew. Math. {\bf 417} (1991), 191-213.

\bibitem[Vis98]{M2} A. Vistoli: \emph{The Chow ring of
$\mathcal{M}_2$}, Invent. Math. {\bf 131} (1998), 635--644, an
appendix to \emph{Equivariant Intersection Theory}, by D. Edidin
and W. Graham.


\end{thebibliography}
\end{document}